\documentclass[leqno]{amsart}
\usepackage{amsmath,amsfonts,amssymb,amscd,amsthm,amsbsy,upref}
\numberwithin{equation}{section}
\newtheorem{thm}{Theorem}[section]
\newtheorem{lem}[thm]{Lemma}
\newtheorem{cor}[thm]{Corollary}
\newtheorem{prop}[thm]{Proposition}
\theoremstyle{definition}
\newtheorem{defn}[thm]{Definition}
\newtheorem*{Defn}{Definition}
\newtheorem*{Defns}{Definitions}
\newtheorem*{remark}{Remark}
\def\ep{\varepsilon}
\def\w{\omega}
\def\bA{\overline{\A}}
\def\tZ{\widetilde Z}
\def\ctA{\widetilde{\A}}
\def\tz{\tilde z}
\def\A{{\mathcal A}}
\def\B{{\mathcal B}}
\def\C{{\mathcal C}}
\def\F{{\mathcal F}}
\def\M{{\mathcal M}}
\def\T{{\mathcal T}}
\def\nat{{\mathbb N}}
\def\real{{\mathbb R}}
\def\Zs{Z^{(*)}}
\def\auc{\operatorname{auc}}
\def\aus{\operatorname{aus}}
\def\cof{\operatorname{cof}}
\def\dist{\operatorname{dist}}
\def\FDD{\operatorname{FDD}}
\def\Sz{S_z}
\def\coo{c_{00}}
\def\supp{\operatorname{supp}}
\def\deltab{\overline{\delta}}
\def\epb{\overline{\ep}}
\def\etab{\overline{\eta}}

\def\hangbox to #1 #2{\vskip1pt\hangindent #1\noindent \hbox to
#1{#2}$\!\!$}
\def\myitem#1{\hangbox to 80pt {#1\hfill}}

\begin{document}
\title{Embedding into Banach spaces with finite dimensional decompositions}
\author{E. Odell \and Th. Schlumprecht}
\address{Department of Mathematics\\ The University of Texas at Austin\\
1 University Station C1200\\ Austin, TX 78712-0257}
\email{odell@math.utexas.edu}
\address{Department of Mathematics\\ Texas A\&M University\\
College Station, TX 77843-3368} \email{schlump@math.tamu.edu}
\thanks{Research supported by the National Science Foundation.}
\begin{abstract}
This paper deals with the following types of problems: Assume
a Banach space $X$ has some property (P). Can it be embedded into
some Banach space $Z$  with a finite dimensional decomposition
having property (P), or more generally, having a property
related to (P)? Secondly, given a class of Banach spaces,
does there exist a Banach space in this class, or in
a closely related one, which is universal for this class?
\end{abstract}
\maketitle

\section{Introduction}\label{S:1}

The fact that every separable infinite dimensional real Banach space
$X$ embeds into $C[0,1]$ dates back to the early days of Banach
space theory \cite[Th\'eor\`eme 9, page 185]{Ba}. This result has
inspired two types of problems. First, given a space $X$ in a
certain class can it be embedded isomorphically into a space $Y$ of
the same class with a basis or, more generally, a finite dimensional
decomposition (FDD)? Secondly, given a class of spaces does there
exist a universal space $X$ for that class which is in the class or
in a closely related one? By saying $X$ is universal for a class $C$
we mean that each $Y\in C$ embeds into $X$. As it happens these two
types of problems are often related in that solving a problem of the
first type can lead to a solution to the analogous problem of second
type.

For example, J.~Bourgain \cite{Bo} asked if there exists a separable
reflexive space $X$ which is universal for the class of all
separable superreflexive Banach spaces. This question arose from his
result that if $X$ contains an isomorph of all separable reflexive
spaces then $X$ is universal, i.e., contains an isomorph of
$C[0,1]$. This improved an earlier result of Szlenk \cite{Sz} who
showed $X^*$ was not separable. Work by S.~Prus \cite{Pr} showed
that it sufficed to prove that for a separable superreflexive space
$Y$ there exists $1<q\le p<\infty$, $C<\infty$ and a space $Z$ with
an FDD $E = (E_i)$ satisfying $C$-$(p,q)$-estimates,
$$C^{-1} \bigg(\sum \|z_i\|^p\bigg)^{1/p}
\le \|\sum z_i\| \le C\bigg(\sum \|z_i\|^q\bigg)^{1/q}$$ for all
block sequences $(z_i)$ of $Z$ w.r.t.\ $(E_i)$. Such a space $Z$ is
automatically reflexive and thus we have the problem of given $p,q$,
when does a reflexive space $Y$ embed into such a space $Z$.

An earlier version of this problem  was raised by
W.B.~Johnson \cite{J1} resulting from his work on $L_p$ and earlier
work with M.~Zippin \cite{JZ1,JZ2}. The problem addressed in
\cite{J1} was to characterize when a subspace $X$ of $L_p$, $1<p<2$,
embeds into $\ell_p$. In \cite{JO} it was shown that if a subspace
$X$ of $L_p$, with $2<p<\infty$, embeds into $\ell_p$ if and only
if  $X$ does not contain an isomorph of $\ell_2$ (later
improved to $X$ almost isometrically embeds into $\ell_p$
\cite{KW}). This characterization does not work in $L_p$, $1\le
 p<2$, since $L_q$ embeds into $L_p$ if $p\le q\le 2$, but the $p>2$
characterization is equivalent (by \cite{KP}) to every
normalized basic sequence in $X$ has a subsequence 2-equivalent to
the unit vector basis of $\ell_p$. Johnson showed that this
criterion (with ``2-equivalent'' replaced  by $C$-equivalent for
some $C<\infty$) characterized when $X\subseteq L_p$, $1<p<2$,
embeds into $\ell_p$. His argument showed that $X$ embedded into
$(\sum H_n)_{\ell_p}$ for some blocking $(H_n)$ of the Haar basis into an
FDD and of course $(\sum H_n)_{\ell_p}$ embeds into
$\ell_p$. Johnson also considered the dual problem which brought
quotient characterizations into the picture. These had appeared
earlier \cite{JZ2} when it was shown that $X$ embeds into $(\sum
E_n)_{\ell_p}$, where $(E_n)$ is a sequence of finite dimensional
Banach spaces iff $X$ is a quotient of such a space.

It turns out that the characterization required to ensure that a
reflexive space $X$ embeds into one with an FDD satisfying
$(p,q)$-estimates is not a subsequence criterion in the general
setting, i.e., if we do not assume $X$ to be a subspace of $L_p$, but
rather one that can be expressed in terms of weakly null trees in
$S_X$, the unit sphere of $X$. This can be viewed as an infinite
version of the notion of asymptotic structure \cite{MMT}. If $X$ is
a Banach space then, for $n\in\nat$, a normalized monotone basis is
 said to be in
the $n^{\text{th}}$-{\em asymptotic structure of} $X$,  and we write
$(e_i)_{i=1}^n \in \{X\}_n$, if for all $\ep>0$ the following holds
($\cof(X)$ will denote the set of all closed subspaces
 of $X$ having finite codimension):
\begin{align}\label{asymptstructure}
&\forall\ X_1\!\in\!\cof(X)\,\exists\ x_1\!\in\!S_{X_1}\ \forall\ 
X_2\!\in\!\cof(X)\,\exists\ x_2\!\in\!S_{X_2}\,\ldots\ \forall\ 
X_n\!\in\!\cof(X)\,\exists\ x_n\!\in\!S_{X_n}\\
&\qquad (x_i)_{i=1}^n\text{ is }(1+\ep)\text{-equivalent to
}(e_i)_{i=1}^n\ .\notag
\end{align}
The fact that  some normalized monotone basis $(e_i)_{i=1}^n$ is a
member of
 $\{X\}_n$ can be, maybe more intuitively, described by a game between two players.
Player I chooses $X_1\in \cof(X)$, then Player II chooses $x_1\in
S_{X_1}$. This procedure
 is repeated until a sequence $(x_i)_{i=1}^n$ is obtained. Player II is declared winner
 of the game if $(x_{i})_{i=1}^n$ is $(1+\ep)$-equivalent to $(e_i)_{i=1}^n$.
Condition \eqref{asymptstructure} means that Player II has a winning
strategy.

 It is not hard to show that
$\{X\}_n$ is a compact subset of $\M_n$, the set of all such
normalized monotone bases $(e_i)_{i=1}^n$ under the metric $\log d_b
(\cdot,\cdot)$ where $d_b((e_i)_{i=1}^n,(f_i)_{i=1}^n)$ is the basis
equivalence constant between the bases.  Lembergs \cite{L} proof of
Krivine's theorem shows
 that there is a $1\le p\le \infty$, so that the unit vector basis of $\ell_p^n$
 is in $\{X\}_n$ for all $n\in\nat$.
In \cite{MMT} it is shown that $\{X\}_n$ is also the smallest closed
subset $\C$ of $\M_n$ with the property that, for all $\ep>0$,
player~I has a winning strategy for forcing player~II to select
$(x_i)_{i=1}^n$ with $d_b((x_i)_{i=1}^n,\C) <1+\ep$. This does not
generalize to produce say $\{X\}_\infty$ since we lose compactness.
However we can still consider a class $\A$ of normalized monotone
bases with the property that in the infinite game player~I has a
winning strategy for forcing II to select $(x_i)_{i=1}^\infty \in
\A$.

These notions can be restated in terms of weakly null trees when
$X^*$ is separable. Indeed $\{X\}_n$ is the smallest class such that
every weakly null tree of length $n$ in $S_X$ admits a branch
$(x_i)_{i=1}^n$ with $d_b((x_i)_{i=1}^n, \{X\}_n) <1+\ep$.
 Precise definitions of weakly null trees and other terminology
appear in Section \ref{S:2}.

If $\A$ is as above for $X$ we can also restate the winning strategy
for player~I in terms of weakly null trees (of infinite level) but
there are some difficulties. First given plays $X_1,X_2,\ldots$ by
player~I we cannot select a branch $(x_i)$ with $x_i\in X_i$ for all
$i$ but only that $x_i$ is close to an element of $S_{X_i}$.
Secondly not all games are determined so we need a fattening
$\A_\ep$ of $\A$ and then need to close it to $\bA_\ep$ in the
product of the discrete topology on $S_X$ to obtain a determined
game. This will lead to the property that if every weakly null tree
in $X$ admits a branch in $\A$ then if $X\subseteq Z$, a space with
an appropriate FDD $(E_i)$, one can find a blocking $(F_i)$ of
$(E_i)$ and $\bar \delta = (\delta)$, $\delta_i\downarrow 0$, so
that every $(x_i)\subseteq S_X$ which is a $\bar\delta$-skipped
block sequence w.r.t.\ $(F_i)$ is in $\bA_\ep$. These will be
defined precisely in Section~2.

An application will be the solution of Johnson's problem (when does
a reflexive space $X$ embed into an $\ell_p$-FDD?)), 
Johnson and Prus' problem
(when does a reflexive space $X$ embed into one with an FDD
satisfying $(p,q)$-estimates) and, as a consequence, Bourgain's
problem. These solutions will be given in Sections~4 and 5. Among
other characterizations we will show that if for some $C<\infty$
every weakly null tree in a reflexive space $X$ admits a branch 
$C$-dominating the
unit vector basis of $\ell_p$ and a branch $C$-dominated by the unit
vector basis of $\ell_q$  then $X$ embeds into a space with an FDD
satisfying $(p,q)$-estimates.

The machinery developed  in Section 2 also has applications in the
nonreflexive setting. In Section~3 we consider and characterize
spaces $X$ of Szlenk index $\w$, the smallest possible. 
If $X$ is a separable Banach space not containing $\ell_1$ then $\Sz(X)$ is
an ordinal index which is less than $\w_1$ iff $X^*$ is separable.
For $\ep >0$ set $K_0 (X,\ep) = B_{X^*}$ and for $\alpha <\w_1$ we
recursively define
\begin{equation*}
K_{\alpha +1} (X,\ep)
 = \left\{ x^* \in K_\alpha (X,\ep) : 
 \begin{matrix} &\exists\ (x_n^*) \subseteq
K_\alpha (x,\ep)\text{ with }\\
\noalign{\vskip4pt}
&w^* - \lim_{n\to\infty} x_n^* = x^*\text{ and }
\liminf_{n\to\infty} \|x_n^* - x^*\| \ge \ep \end{matrix}\right\}
\end{equation*}
If $\alpha$ is a limit ordinal,
$$K_\alpha (X,\ep) = \bigcap_{\beta <\alpha} K_\beta (X,\ep)\ .$$
$S_z(X,\ep)$  is the smallest $\alpha$ with $K_\alpha
(X,\ep)=\emptyset$ or $\w_1$ otherwise.
$$S_z(X) = \sup \{S_z (X,\ep) :  0  < \ep <1\}\ .$$
(This definition is an equivalent version of Szlenk's original index
using Rosenthal's $\ell_1$ theorem \cite{R}.)

We will show that $\Sz(X) =\w$ iff $X^*$ can be embedded as a
$w^*$-closed subspace of a space $Z$ with an FDD satisfying
$1$-$(p,1)$-estimates. A long list of further equivalent conditions
(Theorem~3.4) will be given including that $X$ can be renormed to be
$w^*$-uniform Kadec Klee and $X$ can be renormed to be
asymptotically uniformly smooth (of power type $q$ for some $q>1$).

Asymptotic uniformly smooth (a.u.s.) and asymptotic uniformly convex
(a.u.c) norms, defined in Section~3, are asymptotic versions of
uniformly smooth and uniformly convex due to \cite{JLPS} based upon
modulii of V.D.~Milman \cite{Mi}. Theorem~3.4, mentioned above,
gives the result that $X$ can be given an a.u.s.\ norm iff it can be
given one of power type $q$ for some $q>1$. We obtain a similar
result for a.u.c.\ for reflexive spaces. Recall that Pisier
\cite{Pi} proved that a superreflexive (equivalently, uniformly
convex) space can be renormed to be uniformly convex of power type
$p$ for some $2\le p<\infty$ and similarly for uniformly smooth with
$1<p\le 2$.

In Section 3 we  also give a proof of Kalton's theorem \cite{K} that
a Banach space $X$ embeds into $c_0$ if for some $C<\infty$ every
weakly null tree in $S_X$ admits a branch $(x_i)_{i=1}^\infty$
satisfying $\sup_n \|\sum_1^n x_i\| \le C$. This proof fits nicely
into our Section~2 machinery.

In Section 5 we discuss applications of our results to universal
problems. In regard to Bourgain's problem we show the space
constructed is universal for the class
$$\{ X : X\text{ is reflexive, } \Sz(X) = \Sz (X^*) = \w\}\ ,$$
which includes all superreflexive spaces. We also discuss the
universal problem for reflexive a.u.s.\ (or a.u.c.) spaces.

A central theme of the problems we have presented is
coordinatization. A coordinate-free property is considered and we
wish to embed a space $X$ with this property into a space $Z$ with
an FDD which realizes this property w.r.t.\ its  ``coordinates''.
The tools we use, in addition to the ones mentioned above, are
several. There are the blocking arguments of Johnson and Zippin
\cite{J1}, \cite{JZ1,JZ2} and some known embedding theorems which we
cite now.

\subsection*{1.1 \cite{DFJP}}
If $X^*$ is separable then $X$ is a quotient of a space with a
shrinking basis.

\subsection*{1.2 \cite{Z}}
If $X^*$ is separable then $X$ embeds into a space with a shrinking
basis.

\subsection*{1.3 \cite{Z}}
If $X$ is reflexive then $X$ embeds into a reflexive space with a
basis.

We will often begin with $X\subseteq Z$, one of the spaces given by
1.2, 1.3 or with $X$ a quotient of $Z$ (as in 1.1) and the problem
will be to put a new norm on $Z$ which reflects the structure of $X$
that we wish to coordinatize and maintains that $X$ is a subspace of
$Z$ (or a quotient).

All of our Banach spaces in this paper are real  and separable. We
will use $X,Y,Z,\ldots$ for infinite dimensional spaces and
$E,F,G,\ldots$ for finite dimensional spaces or write $E= (E_n)$ 
for an FDD.

Most of the results we will present have appeared in a number of
recent papers  (\cite{OS1}, \cite{OS2}, \cite{OSZ} \cite{KOS}, \cite{K},
\cite{GKL}, \cite{JLPS}). As the theory has developed the proofs and
results have been better understood, generalized and improved. Our
aim is to give a unified presentation of these improvements and
in several cases present easier proofs. New results are also included.


\section{A general combinatorial result}\label{S:2}

In this section we state and prove three general combinatorial
results (Theorem~\ref{T:2.2} and Corollaries \ref{C:2.6} and \ref{C:2.8}). 
These are reformulations and improvements of results in \cite{OS1}.
We will present a different more accessible proof.

We first introduce some notation.

Let $Z$ be a Banach space with an FDD $E=(E_n)$.  For $n\in\nat$ we
denote the $n$- th {\em coordinate projection } by $P^E_n$, i.e.
$P_n^E:Z\to E_n,\quad \sum z_i\mapsto z_n.$ For  finite  $A\subset
\nat$ we put $P^E_A=\sum_{n\in A} P_n^E$. The  {\em projection constant  }
of $(E_n)$ (in $Z$) is defined by
$$K=K(E,Z)=\sup_{m\le n}\|P_{[m,n]}^E\|\ .$$

Recall that $K$ is always finite and, as in the case of bases, we call
{\em  $(E_n)$ bimonotone (in $Z$)}
if $K=1$.  By passing to the equivalent norm
$$|||\cdot|||:Z\to\real,\quad z\mapsto\sup_{m\le n} \|P^E_{[m,n]}(z)\|\ ,$$
we can always renorm $Z$ so that $K=1$.

For a sequence $(E_i)$ of finite dimensional spaces we define the
vector space
$$\coo(\oplus_{i=1}^\infty E_i)=
\Big\{(z_i):z_i\in E_i, \text{ for $i\in\nat,$ and $\{i\in\nat:
z_i\not=0\}$ is finite}\big\}\ ,$$ 
which is dense in each Banach space for which $(E_n)$ is an FDD. 
For $A\subset \nat$ we denote  by $\oplus_{i\in A} E_i$ the linear subspace
of $\coo(\oplus E_i)$ generated by the elements
of $(E_i)_{i\in A}$ and we denote its closure in $Z$ by $(\oplus E_i)_Z$.
As usual we denote the vector space of sequences
in $\real$ which are eventually zero by $c_{00}$ and its  unit 
vector basis by $(e_i)$.

The vector space $\coo(\oplus_{i=1}^\infty E^*_i)$,
where $E^*_i$ is the dual space of $E_i$, for $i\in\nat$,  is  a $w^*$-dense
subspace of $Z^*$. 
(More precisely $E_i^*$ is the subspace of $Z^*$ generated by all elements
$z^*$ for which $z^*|_{E_n} =0$ if $n\ne i$. 
$E_i^*$ is uniformly isomorphic to the dual space of $E_i$ and is 
isometric to it if $K(E,Z)=1$.) 
We denote the norm closure of $\coo(\oplus_{i=1}^\infty E^*_i)$
in $Z^*$ by $\Zs$. 
$\Zs$ is $w^*$-dense in $Z^*$, the unit ball $B_{\Zs}$ norms $Z$  and 
$(E_i^*)$ is an FDD of $\Zs$ having a projection constant not 
exceeding $K(E,Z)$. 
If $K(E,Z)=1$ then $B_{Z^{(*)}}$ is 1-norming and $Z^{(*)(*)}=Z$.

For $z\in \coo(\oplus E_i)$ we define the $E$-{\em support of $z$} by
$$\supp_E(z)=\big\{i\in\nat:  P^E_i(z)\not=0\big\}\ .$$
A  non-zero sequence (finite or infinite) $(z_j)\subset \coo(\oplus E_i)$  is
called  a {\em block sequence of $(E_i)$} if
$$\max\supp_E(z_n)<  \min\supp_E(z_{n+1}),\text{ whenever $n\in\nat$ (or
$n\!<\!\text{length}(z_j)$)},$$
and it is called a  {\em skipped block sequence of $(E_i)$} if
$1<\min\supp_E(z_1)$ and 
$$\max\supp_E(z_n) <  \min\supp_E(z_{n+1})-1,
\text{ whenever $n\in\nat$ (or $n<\text{length}(z_i)$)}.$$
Let $\deltab=(\delta_n)\subset (0,1]$.
A (finite or infinite) sequence $(z_j)\subset S_Z=\{z\in Z:
\|z\|=1\}$ is called a $\deltab$-{\em block sequence of $(E_n)$}
or a $\deltab$-{\em skipped block sequence of $(E_n)$}
if there are $1\le k_1< \ell_1< k_2< \ell_2<\cdots$ in $\nat$ so that
$$\| z_n - P^E_{[k_{n},\ell_n]}(z_n)\|<\delta_n,\text{ or }
\| z_n - P^E_{(k_{n},\ell_n]}(z_n)\|<\delta_n, \text{ respectively,}$$
for all $n\in\nat$  (or $n\le\text{length}(z_j)$).
Of course one could generalize the notion of $\deltab$-block and
$\deltab$-skipped block sequences to more general sequences, but
we prefer to introduce this notion only for normalized sequences.
It is important to note that in the definition of $\deltab$-skipped
block sequences  $k_1\ge 1$, and that therefore the $E_1$-coordinate of
$z_1$ is small (depending on $\delta_1$).

A sequence of finite-dimensional spaces $(G_n)$ is  called
a {\em blocking of $(E_n)$} if there are
$0=k_0<k_1<k_2<\cdots$ in $\nat$ so that
$G_n=\oplus_{i=k_{n-1}+1}^{k_n} E_i$, for $n=1,2,\ldots$.

We denote the sequences in $S_Z$ of length $n\in\nat$ by $S_Z^n$ and
the infinite sequences in $S_Z$  by $S^\omega_Z$. 
For $m,n\in\nat$, for $x=(x_1,x_2,\ldots, x_m)\in S^m_Z$
and $y=(y_1,y_2,\ldots, y_n)\in S^n_Z$ or $y=(y_i)\in S^\omega_Z$ we denote
the concatenation of $x$ and $y$  by $(x,y)$, i.e.
$$(x,y)=(x_1,x_2,\ldots, x_m,y_1,\ldots, y_m),\text{ or }
(x,y)=(x_1,x_2,\ldots, x_m,y_1,y_2,\ldots) \text{ respectively }.$$
We also allow  the case $x=\emptyset$ or $y=\emptyset$  and
let $(\emptyset,y)=y$ and $(x,\emptyset)=x$.

Let $\A\subset S_Z^\w$ be given. 
We denote the closure of $\A$ with respect to the product topology of 
the discrete topology on $S_Z$ by $\overline{\A}$. 
Note that if $\A$ is closed it follows that for $x=(x_i)\in S_Z^w$,
\begin{equation}\label{E:2.1}
x\in\A\iff \forall\ n\in\nat\ \exists\ z\in S_Z^\omega\quad
(x_1,x_2,\ldots, x_n, z)\in\A
\end{equation}

If  $\epb=(\ep_i)$ is a sequence in $[0,\infty)$ we write
$$\A_{\epb}=\big\{(z_i)\in S^\omega_Z:\ \exists\  (\tz_i)\in\A,\quad 
\|z_i-\tz_i\|\le \ep_i\text{ for all } i\in\nat\big\}$$
and call the set $\A_{\epb}$ {\em the $\epb$-fattening of $\A$}.
For $\ell\in\nat$ and  $\epb=(\ep_i)_{i=1}^\ell\subset[0,\infty)$ 
we let $\A_{\epb}=\A_{\deltab}$, where $\deltab=(\delta_i)$ and
$\delta_i=\ep_i$, for $i=1,2,\ldots, \ell$ and $\delta_i=0$ if $i>\ell$.

If  $\ell\in\nat$ and  $x_1,x_2,\ldots, x_\ell\in S_Z$ we let
$$\A(x_1,x_2,\ldots, x_\ell)=
\big\{z=(z_i)\in S_Z^\omega: (x_1,x_2,\ldots, x_\ell,z)\in\A\big\}\ .$$

Let $\A\subset S^\omega_Z$ and
$\B=\prod_{i=1}^\infty B_i$, where $B_n\subset S_Z$ for $n\in\nat$.

We consider the following {\em $(\A,\B)$-game } between two players: 
Assume that $E=(E_i)$ is an FDD for $Z$.
\begin{align*}
&\text{Player I chooses }n_1\in\nat\ ,\\
&\text{Player II chooses } z_1\in \coo\big(\oplus_{i=n_1+1}^\infty
E_i\big)\cap B_1\ ,\\
&\text{Player I chooses }n_2\in\nat\ ,\\
&\text{Player II chooses } z_2\in\coo\big(\oplus_{i=n_2+1}^\infty
E_i\big)\cap B_2\ ,\\
&\vdots
\end{align*}
Player I wins  the $(\A,\B)$-game if the resulting sequence $(z_n)$
lies in $\A$. 
If Player I has a winning strategy
(forcing the sequence $(z_i)$ to be in $\A$) we will write
$WI(\A,\B)$ and if Player II has a winning strategy
(being able to choose $(z_i)$ outside of $\A$) we write $WII(\A,\B)$.  
If $\A$ is a Borel set with respect
to the product  of the discrete topology on $S^\omega_Z$
(note that $\B$ is always closed in
the product  of the discrete topology on $S^\omega_Z$), it follows
from the main   theorem in \cite{Ma} that the game is determined, i.e., either
$WI(\A,\B)$ or $WII(\A,\B)$.

Let us define $ WII(\A,\B)$ formally.
We will need to introduce trees in Banach spaces.

We define
\begin{equation*}
T_\infty=\bigcup_{\ell\in\nat}\big\{(n_1,n_2,\ldots,n_\ell)
 :n_1<n_2<\cdots n_\ell \text{ are  in $\nat$}\big\}\ .
\end{equation*}
If $\alpha=(m_1,m_2,\ldots, m_\ell)\in T_\infty$, we call $\ell$ the
{\em  length of $\alpha$} and denote it by $|\alpha|$, and 
$\beta=(n_1,n_2,\ldots, n_k)\in T_\infty$ is called
an {\em extension of $\alpha$}, or $\alpha$ is called {\em a restriction of} 
$\beta$, if $k \ge \ell$ and $n_i=m_i$, for $i=1,2,\ldots,\ell$.
We then write $\alpha\le \beta$ and with this order $(T_\infty,\le)$ is a tree.

In this work  {\em trees} in a Banach space $X$ are families  in
$X$ indexed by $T_\infty$,  thus they are countable infinitely branching trees
of countably infinite length.

For a  tree $(x_\alpha)_{\alpha\in T_\infty}$ in a  Banach space $X$,
and $\alpha=(n_1,n_2,\ldots, n_\ell)\in T_\infty\cup\{\emptyset\}$
we call the sequences of the form $(x_{(\alpha,n)})_{n>n_\ell}$
{\em nodes of $(x_\alpha)_{\alpha\in  T_\infty}$}.
The sequences $(y_n)$,
with $y_i=x_{(n_1,n_2,\ldots,n_i)}$,  for $i\in\nat$,  for some
strictly increasing sequence $(n_i)\subset\nat$, are called {\em branches of }
$(x_\alpha)_{\alpha\in T_\infty}$. 
Thus, branches of a tree $(x_\alpha)_{\alpha\in T_\infty}$ are sequences of 
the form $(x_{\alpha_n})$ where $(\alpha_n)$ is a maximal linearly
ordered (with respect to extension) subset of $T_\infty$.

If $(x_\alpha)_{\alpha\in T_\infty}$  is a tree in $X$ and if 
$T'\subset T_\infty$  is closed under
taking restrictions so that for each $\alpha\in T'\cup\{\emptyset\}$
infinitely many direct successors of $\alpha$ are also in $T'$ then we
call $(x_\alpha)_{\alpha\in T'}$ a {\em full subtree of}
$(x_\alpha)_{\alpha\in T_\infty}$. 
Note that
$(x_\alpha)_{\alpha\in T'}$ could then be relabeled to a family indexed by
$T_\infty$ and note that the branches of $(x_\alpha)_{\alpha\in T'}$
are branches of $(x_\alpha)_{\alpha\in T_\infty}$ and that the nodes of 
$(x_\alpha)_{\alpha\in T'}$ are subsequences of  certain nodes of
$(x_\alpha)_{\alpha\in T_\infty}$.

We call a tree $(x_\alpha)_{\alpha\in T_\infty}$  in a  Banach space $X$
{\em normalized} if $\|x_\alpha\|=1$, for all  $\alpha\in T_\infty$ and
{\em weakly null} if every node is weakly null. 
More generally if $\T$ is a topology on $X$ and
a tree  $(x_\alpha)_{\alpha\in T_\infty}$  in a  Banach space $X$
is called {\em $\T$-null } if every node converges to $0$ with respect to $\T$.

If  $(x_\alpha)_{\alpha\in T_\infty}$ is a tree in a Banach space $Z$
which has an FDD $(E_n)$  we call it a
{\em  block tree  of } $(E_n)$ if every node is a block sequence of $(E_n)$.

We will also need to consider trees of finite length.  
For $\ell\in\nat$ we
call a family $(x_\alpha)_{\alpha\in T_\infty,|\alpha|\le \ell}$ in
$X$ a {\em tree of length $\ell$}. 
Note that the notions nodes, branches, $\T$-null and block trees 
can be defined analogously for trees of finite length.

\begin{Defn} 
Assume that $Z$ is a Banach space with
an FDD $(E_i)$, $\A\subset S_Z^\omega$ and $\B=\prod_{i=1}^\infty B_i$, 
with $B_i\subset S_Z$ for $i\in\nat$. 
We say that
{\em Player II has a winning strategy for the $(\A,\B)$-game} if

\myitem{($WII(\A,\B)$)} 
There exists a  block tree $(x_\alpha)_{\alpha\in T_\infty}$  of $(E_i)$ in 
$S_Z$ all of whose branches are in $\B$ but none of its branches are in $\A$.

In case that the $(\A,\B)$-game is determined $WI(\A,\B)$  can be 
therefore stated as follows.

\myitem{($WI(\A,\B)$)} 
Every block tree $(x_\alpha)_{\alpha\in T_\infty}$
of $(E_i)$ in $S_X$, all of whose branches are in $\B$, has a branch in $\A$.
\end{Defn}

The proof of the following Proposition is easy.

\begin{prop}\label{P:2.1}
Let $\A,\ctA \subset S^\omega_Z$, $\B=\prod_{i=1}^\infty B_i$,
with $B_i\subset S_Z$ for $i\in\nat$. Assume that the
$(\A,\B)$-game and the $(\ctA,\B)$-game are determined.
\begin{enumerate}
\item[a)] If $\A\subset \ctA$, then
  $$WI(\A,\B)\Rightarrow WI(\ctA,\B)\text{ and }
WII(\ctA,\B)\Rightarrow WII(\A,\B)\ .$$
\item[b)]  $WI(\A,\B)\iff \exists\  n\!\in\nat \ \forall\ 
x\!\in\!\big(\oplus_{i=n+1}^\infty E_i\big)\cap B_1\quad
WI(\A(x),\prod_{i=2}^\infty B_i)$
 \item[c)] If $\ell\in\nat$,
$\epb=(\ep_i)_{i=1}^\ell\subset[0,\infty)$ and $x_i,y_i\in B_i$
 with  $\|x_i-y_i\|\le \ep_i$ for $i=1,2,\ldots, \ell$ then
$$ WI\Big(\A(x_1,x_2,\ldots,x_\ell),\prod_{i=\ell+1}^\infty B_i\Big)\Rightarrow
WI\Big(\A_{\epb}(y_1,y_2,\ldots,y_\ell), \prod_{i=\ell+1}^\infty B_i\Big)\ .$$
\end{enumerate}
\end{prop}

\begin{lem}\label{L:2.1a} 
Let $\A$, and $\epb=(\ep_i),\deltab=(\delta_i)\subset [0,\infty)$
Then
$$\overline{\big(\,\overline{\A_{\epb}}\big)_{\deltab}}
\subset\overline{\A_{\epb+\deltab}}\ .$$
\end{lem}

\begin{proof} We observe
\begin{align*}
u=(u_i)&\in \overline{\big(\,\overline{\A_{\epb}}\big)_{\deltab}} \\
\iff& \forall\ n\in\nat\ \exists\  v^{(n)}\in S_Z^\omega 
\quad (u_1,\ldots,u_n,v^{(n)})\in \big(\,\overline{\A_{\epb}}\big)_{\deltab}\\
\Longrightarrow\ &\forall\ n\in\nat\ \exists\ x_1,x_2,\ldots, x_n\in S_Z 
\text{ and } w^{(n)}\in S_Z^\omega\\
&\|x_i-u_i\|\le \delta_i, \text{ for }i=1,\ldots, n,\text{ and } 
(x_1,\ldots,x_n,w^{(n)})\in \overline{\A_{\epb}}
\\
\Longrightarrow\ &\forall\ n\in\nat\ \exists\  x_1,x_2,\ldots, x_n\in S_Z
\text{ and } w^{(n)}\in S_Z^\omega
\  \forall\ m\in\nat\  \exists\  y^{(m)}\in S_Z^\omega\\
\qquad&\|x_i-u_i\|\le \delta_i\text{ for }i=1,2,\ldots, n \text{ and }
(x_1,\ldots,x_n,w_1^{(n)},w_2^{(n)},\ldots,
w_m^{(n)},y^{(m)})\!\in\! \A_{\epb}\\
\Longrightarrow\ &\forall\ n\in\nat\  \exists\  x_1,x_2,\ldots,
x_n \in S_Z\  \exists\  y^{(n)}\in S_Z^\omega\\
\qquad&\|x_i-u_i\|\le \delta_i\text{ for }i=1,2,\ldots, n \text{ and }
(x_1,\ldots,x_n,y^{(n)})\!\in\! \A_{\epb}\\
\Longrightarrow\ &\forall\ \ell\!\in\!\nat\ 
\exists\ z^{(\ell)}\in \A \quad \|u_i-z_i^{(\ell)}\|\le\delta_i+\ep_i, 
\text{ for }i=1,2,\ldots ,\ell\\
 \iff  &u\in \overline{\A_{\epb+\deltab}}\ .
\end{align*}\end{proof}

Now we can state one of our main combinatorial principles.

\begin{thm}\label{T:2.2} 
Let $Z$ have an FDD $(E_i)$ and let $B_i\subset S_Z$, for $i=1,2,\ldots $. 
Put $\B=\prod_{i=1}^\infty B_i$ and let $\A\subset S^\omega_Z$.

Assume that for all $\epb=(\ep_i)\subset(0,1]$ we have
$WI(\,\overline{\A_{\epb}}, \B)$.

Then  for all  $\epb=(\ep_i)\subset(0,1]$ there exists a
blocking $(G_i)$ of $(E_i)$ so that every skipped block sequence $(z_i)$ of
$(G_i)$, with $z_i\in B_i$, for $i\in\nat$, is in $\overline{\A_{\epb}}$.
\end{thm}

\begin{proof}
Let $\epb=(\ep_i)\subset(0,1]$ be given. For $k=0,1,2,\ldots$  put
$\epb^{(k)}=(\ep_i^{(k)})$ with \\  $\ep_i^{(k)}=\ep_i(1-2^{-k})/2$
for $i\in\nat$.

We put $\ctA=\overline{\A_{\epb/2}}$.

For $\ell\in\nat$ we write
$\B^{(\ell)}=\prod_{i=\ell+1}^\infty B_i$.

By induction we choose for $k\in\nat$ numbers $n_k\in\nat$ so that
$0=n_0<n_1<n_2<\cdots$, and so that for any $k\in\nat$, if 
$G_k=\oplus_{i=n_{k-1}+1}^{n_k} E_i$,
\begin{align}\label{E:2.2.1}
&WI\big(\,\overline{\ctA}_{\epb^{(k)}}(\sigma,x),\B^{(\ell+1)}\big)
\text{ for any $0\le\ell<k$ and  any normalized skipped block}\\
&\text{
 $\sigma=(x_1,x_2,\ldots ,x_\ell)\in\prod_{i=1}^\ell B_i$ of $(G_i)_{i=1}^{k-1}$
 ($\sigma=\emptyset$ if $\ell=0$)}\notag\\
 &\text{and $x\in S_{\oplus_{i=n_k+1}^\infty E_i}\cap B_{\ell+1}$}\notag\\
\label{E:2.2.2} &WI
\big(\,\overline{\ctA}_{\epb^{(k)}}(\sigma),\B^{(\ell)}\big) \text{
for any $0\le\ell<k$ and any
normalized skipped block} \\
&\text{$\sigma=(x_1,x_2,\ldots, x_\ell)\in \prod_{i=1}^\ell B_i $ of
$(G_i)_{i=1}^{k}$}\notag
\end{align}
For $k=1$ we deduce from Proposition~\ref{P:2.1} (b), Lemma~\ref{L:2.1a} 
and the hypothesis that there is an $n_1\in\nat$
so that $WI\big(\,\overline{\ctA}_{\epb^{(1)}}(x),\B^{(1)}\big)$ for any
$x\in  S_{\oplus_{i=n_1+1}^\infty E_i}\cap B_i$.
This implies \eqref{E:2.2.1} and   \eqref{E:2.2.2}
(note that for $k=1$ $\sigma$ can only be chosen to be $\emptyset$
in  \eqref{E:2.2.1} and   \eqref{E:2.2.2}).

Assume $n_1<n_2<\cdots n_k$ have been chosen for some $k\in\nat$. 
We will first choose $n_{k+1}$ so that \eqref{E:2.2.1} is satisfied. 
In the case that $k=1$ we simply choose $n_2=n_1+1$ and note that
\eqref{E:2.2.1}  for $k=2$ follows from \eqref{E:2.2.1}  for $k=1$ 
since in both cases $\sigma=\emptyset$ is the only choice. 
If $k>1$ we can use the compactness of the sphere of a finite dimensional
space and choose a  finite set $\F$ of normalized  skipped blocks
$(x_1,x_2,\ldots, x_\ell)\in \prod_{i=1}^\ell B_i$,
of $(G_i)_{i=1}^k$ so that for any  $\ell\le k$ and any normalized skipped
block  with length $\ell$,
$\sigma=(x_1,x_2,\ldots,x_\ell)\in\prod_{i=1}^\ell B_i $ of
$(G_i)_{i=1}^k$,  there is a $\sigma'=(x_1',x_2',\ldots,x_\ell')\in \F$ with
$\|x_i-x_i'\|<\ep_i 2^{-k-2}$, for $i=1,2,\ldots,\ell$. 
Then, using the induction hypothesis \eqref{E:2.2.2} for $k$,
and Proposition~\ref{P:2.1}~(b), we choose
$n_{k+1}\in\nat$ large enough so that 
$WI\big(\,\overline{\ctA}_{\epb^{(k)}}(\sigma,x), \B^{(\ell+1)}\big)$
for any $\sigma\in \F$  and 
$x\in S_{\oplus_{i=n_{k+1}+1}^\infty E_i}\cap B_{\ell+1}$.
{From} Proposition~\ref{P:2.1}~(c) and our choice of $\F$ we deduce
$WI\big(\,\overline{\ctA}_{\epb^{(k+1)}}(\sigma,x),\B^{(\ell+1)}\big)$
for any $0\le\ell<k$, any  normalized skipped block $\sigma$ of 
$(G_i)_{i=1}^k$  of length $\ell$ in
$\prod_{i=1}^\ell B_i$ and any 
$x\in S_{\oplus_{i=n_{k+1}+1}^\infty E_i}\cap B_{\ell+1}$,
and, thus, (using the induction hypothesis for $\sigma=\emptyset$) we deduce
\eqref{E:2.2.1}  for $k+1$.

In order to verify    \eqref{E:2.2.2}   let $\sigma=(x_1,x_2,\ldots,
x_\ell)\in \prod_{i=1}^\ell B_i$ be
a normalized skipped block of $(G_i)_{i=1}^{k+1}$ 
(the case $\sigma=\emptyset$ follows from the induction hypothesis). 
Then $\sigma'=(x_1,x_2,\ldots, x_{\ell-1})$ is empty or
a normalized skipped block sequence of $(G_i)_{i=1}^{k-1}$
in $\prod_{i=1}^{\ell-1} B_i$. 
In the second case 
$WI\big(\,\overline{\ctA}_{\epb^{(k+1)}})(\sigma),\B^{(\ell)}\big)=
WI\big(\,\overline{\ctA}_{\epb^{(k+1)}})(\sigma',x_\ell),\B^{(\ell)}\big)$
follows from  \eqref{E:2.2.1} for $k$
and from Proposition~\ref{P:2.1}~(a).
This finishes the  recursive definition of the $n_k$'s and $G_k$'s.

Let $(z_n)$ any normalized skipped block sequence of $(G_i)$
which lies in $\B$. 
For any $n\in\nat$ it follows from \eqref{E:2.2.2} for  $\sigma=(z_i)_{i=1}^n$ 
that  $WI(\,\overline{\ctA_{\epb/2}}(\sigma),\B)$, and, thus,
$\ctA_{\epb/2}(\sigma)\not=\emptyset$, which means
that $\sigma$ is extendable to a sequence in $\ctA_{\epb/2}$
(note that $\lim_{n\to\infty} \ep^{(n)}_i=\ep_i$). 
Thus, any normalized skipped block sequence which is element of $\B$
lies in $\overline{\ctA_{\epb/2}}$ and, thus, by Lemma~\ref{L:2.1a},
in  $\overline{\A_{\epb}}$.
\end{proof}

Now let $X$ be a closed subspace  of $Z$  having an FDD $(E_i)$ and
$\A\subset S_X^\omega$. We consider the following game
\begin{align*}
&\text{Player I chooses }n_1\in\nat\ , \\
&\text{Player II chooses } x_1\in\big(\oplus_{i=n_1+1}^\infty
E_i\big)_Z\cap X,\,\, \|x_1\|=1\ ,\\
&\text{Player I chooses }n_2\in\nat\ ,\\
&\text{Player II chooses } x_2\in\big(\oplus_{i=n_2+1}^\infty
E_i\big)_Z\cap X,\,\, \|x_2\|=1\ ,\\
&\vdots
\end{align*}
As before, Player I wins  if $(x_i)\in \A$. 
Since the game does
not only depend on $\A$ but on the superspace $Z$ in which $X$ is 
embedded and its FDD $(E_i)$ we call this the $(\A,Z)$-game.

\begin{Defn}  
Assume that $X$ is a subspace of a space $Z$ which has an FDD $(E_i)$
and that $\A\subset S^\omega_X$. Define for $n\in\nat$
$$X_n=X\cap \big(\oplus_{i=n+1}^\infty E_i\big)_Z=
\{ x\in X: \forall\ z^*\in\oplus_{i=1}^n E_i^*  \quad z^*(x)=0\}\big)\ , $$ 
a closed subspace of finite codimension in $X$.

We say that {\em Player II has a winning strategy in the
$(\A,Z)$-game if}

\myitem{$WII(\A,Z)$} there is a tree $(x_\alpha)_{\alpha\in
T_\infty}$ in $S_X$ so that
for any $\alpha=(n_1,\ldots, n_\ell)\in T_\infty\cup{\emptyset}$
$x_{(\alpha,n)}\in X_n \text{ whenever }n>n_\ell $, and so that no
branch lies in $\A$.

\noindent In the case that the $(\A,Z)$-game is determined {\em Player I
has a winning strategy in the $(\A,Z)$-game if} the
negation of $WII(\A,Z)$ is true and thus

\myitem{$WI(\A,Z)$} for any tree
 $(x_\alpha)_{\alpha\in T_\infty}$ in $S_X$ so that
 for any $\alpha=(n_1,\ldots, n_\ell)\in T_\infty\cup{\emptyset}$
$x_{(\alpha,n)}\in X_n$ whenever $n>n_\ell $, there is branch in $\A$.
\end{Defn}

For $\A\subset S_X^\omega\subset S_Z^\omega$  and a sequence
$\epb=(\ep_i)$ in $[0,\infty)$ we understand by
$\A_{\epb}$ the $\epb$-fattening of $\A$ as a subset of $S_Z^\omega$. 
In case we want to restrict ourselves to $S_X$ we write $\A^X_{\epb}$, i.e.
$$\A^X_{\epb}=\A_{\epb}\cap S_X^\omega=
\big\{(x_i)\in S_X^\omega:\exists\  (y_i)\!\in\!\A\quad \|x_i-y_i\|\le
\ep_i\text{ for all }i\in\nat\big\}\ .$$
Since $S_X^\omega$ is closed in $S_Z^\omega$ with respect to the
product of the discrete topology, we deduce that
$\overline{\A}^X=\overline{\A^X}$ for $\A\subset S^\omega_X$.

The following Proposition reduces the $(\A,Z)$-game to a game we
treated before. 
In order to be able to do so we need some
technical assumption on the embedding of $X$ into $Z$
(see condition \eqref{E:2.4.1} below).

\begin{prop}\label{P:2.4}  
Let $X\subset Z$, a space with an FDD $(E_i)$. 
Assume the following condition on $X$, $Z$
and the embedding of $X$ into $Z$ is satisfied:
\begin{align} \label{E:2.4.1} 
&\text{ There is a $C>0$ so that for all $m\in\nat$
and $\ep>0$ there is an $n=n(\ep,m)$}\ge m\\
&\qquad\|x\|_{X/X_m}\le C \big[\|P^E_{[1,n]}(x)\|+\ep\big]\text{
whenever }x\in S_X\ .\notag
\end{align}
Assume that $\A\subset S_X^\omega$ and that
for all null sequences  $\epb\subset(0,1]$
we have $WI(\,\overline{\A_{\epb}^X},Z)$.

Then it follows for all null sequences $\epb=(\ep_i)\subset (0,1]$
that $WI( \overline{\A_{\epb}},(S^\omega_X)_{\deltab})$ holds, where
$\deltab=(\delta_i)$ with $\delta_i=\ep_i/28CK$ for
$i\in\nat$,  with $C$ satisfying \eqref{E:2.4.1} and
$K$ being the projection constant of $(E_i)$ in $Z$.
\end{prop}

\begin{proof}
Let $\A\subset S_X^\omega$ and assume  that
$WI(\,\overline{\A^X_{\etab}},Z)$ is satisfied for all null sequences
$\etab=(\eta_i)\subset (0,1]$. 
For a null sequence $\epb=(\ep_i)\subset (0,1]$ we need to verify
$WI( \overline{\A_{\epb}},(S^\omega_X)_{\deltab})$
(with $\delta_i=\ep_i/28KC$ for $i\in\nat$) and so we
let $(z_\alpha)_{\alpha\in T_\infty}$ be a block tree of
$(E_i)$ in $S_Z$ all of whose branches lie in
$(S_X^\omega)_{\deltab}=\{(z_i)\in S^\omega_Z: 
\dist(z_i,S_X)\le \delta_i\text{ for } i=1,2,\ldots\}$.

After passing to a full subtree of $(z_\alpha)$ we can assume
that for any $\alpha=(m_1,\ldots ,m_\ell)$ in $T_\infty$
\begin{equation}\label{E:2.4.2}
z_\alpha\in\oplus_{j=1+n(\delta_\ell,m_\ell)}^\infty E_j
\end{equation}
(where $ n(\ep,m)$ is chosen as in \eqref{E:2.4.1}).

For $\alpha=(m_1,m_2,\ldots, m_\ell)\in T_\infty$ we choose
$y_\alpha\in S_X$ with $\|y_\alpha-z_\alpha\|<2\delta_\ell$
 and, thus, by \eqref{E:2.4.2}
$$\|P^{E}_{[1,n(\delta_\ell,m_\ell)]}(y_\alpha)\|
=\|P^{E}_{[1,n(\delta_\ell,m_\ell)]}(y_\alpha-z_\alpha)\|\le
 2K\delta_\ell\ .$$
Using \eqref{E:2.4.1} we can therefore choose
an $x'_\alpha\in X_{m_\ell}$ so that
$$\|x'_\alpha-y_\alpha\|\le C(2K\delta_\ell+\delta_\ell)\le
3CK\delta_\ell\ ,$$
and thus
$$ 1-3CK\delta_\ell\le \|x'_\alpha\|\le 1+3CK\delta_\ell\ .$$
Letting $x_\alpha=x'_\alpha/\|x'_\alpha\|$ we deduce that
\begin{align*}
\|y_\alpha-x_\alpha\|
&\le \|y_\alpha-x'_\alpha\|+\|x'_\alpha-x_\alpha\|\\
&\le 3CK\delta_\ell + (1+3CK\delta_\ell)3CK\delta_\ell/(1-3CK\delta_\ell)\le
12CK\delta_\ell\notag
\end{align*}
(the last inequality follows from the fact that
$(1+3CK\delta_\ell)/(1-3CK\delta_\ell)\le 3$) and, thus,
$$\|z_\alpha-x_\alpha\|\le 14CK\delta_\ell=\ep_\ell/2\ .$$
Using $WI(\,\overline{\A^X_{\epb/2}}, Z)$ and noting 
that $x_\alpha\in X_{m_\ell}$, for
$\alpha=(m_1,m_2,\ldots, m_\ell)\in T_\infty$ we can choose
a branch of $(x_\alpha)$ which is in $\overline{\A^X_{\epb/2}}$.
Thus, the corresponding branch of $(z_\alpha)$ lies in $\overline{\A_{\epb}}$.
\end{proof}

{From} \cite[Lemma 3.1]{OS1}  it follows that
every separable Banach space $X$ is a subspace of a space $Z$ with an FDD
satisfying the condition \eqref{E:2.4.1} (with $n(m)=m$). 
The following Proposition exhibits two general situations in which   
\eqref{E:2.4.1}  is automatically satisfied.

\begin{prop}\label{P:2.5}
Assume $X$ is a subspace of a space $Z$ having an FDD $(E_i)$. 
In the following two cases \eqref{E:2.4.1} holds:
\begin{enumerate}
\item[a)] If $(E_i)$ is a shrinking FDD for $Z$.
In that case $C$ in \eqref{E:2.4.1} can be chosen arbitrarily close to $1$.
\item[b)] If  $(E_i)$ is boundedly complete for $Z$ (i.e., $Z$ is the dual 
of $\Zs$) and the ball of $X$ is a w$^*$-closed subset of $Z$.
In that case $C$ can be chosen to be the projection constant $K$ of
$(E_i)$ in $Z$.
\end{enumerate}
\end{prop}

\begin{proof}
In order to prove (a) we  will  show that for any  $m\in\nat$ and
any $0<\ep<1$ there is an $n= n(\ep,m)$ so that
$$\|x\|_{X/X_m}\le (1+\ep) \big[\|P^E_{[1,n]}(x)\|+\ep\big],
\text{ whenever }x\in S_X$$
(i.e., $C$  in \eqref{E:2.4.1} can be chosen arbitrarily close to $1$).

Since $X/X_m$ is finite dimensional and
$$(X/X_m)^*=X_m^\perp=\big\{x^*\in X^*: x^*|_{X_m}\equiv0\big\}\ ,$$
we can choose a finite set $A_m\subset S_{X_m^\perp}\subset S_{X^*}$
for which
$$\|x\|_{X/X_m}\le  (1+\ep) \max_{f\in A_m} |f(x)| \text{ whenever }x\in X\ .$$
By the Theorem of Hahn Banach we can extend each $f\in A_m$ to an
element $g\in S_{Z^*}$.
Let us denote the set of all of these extensions $B_m$. 
Since $B_m$ is finite and since  $(E^*_i)$ is an FDD of $Z^*$ we can choose an
$n=n(\ep,m)$ so that $\|P_{[1,n(m)]}^{E^*}(g)-g\|<\ep$ for all $g\in B_m$.
Since $P_{[1,n(m)]}^{E^*}$ is the adjoint operator of
$P_{[1,n(m)]}^{E}$  (consider $P_{[1,n(m)]}^{E^*}$ to be an operator
from $Z^*$ to $Z^*$ and $P_{[1,n(m)]}^{E}$ to be an operator from
$Z$ to $Z$), it follows  for $x\in S_X$,   that
\begin{align*}
\|x\|_{X/X_m}&\le  (1+\ep)\max_{g\in B_m} |g(x)|\\
&\le (1+\ep)\max_{g\in B_m}\big[ \big|P_{[1,n(m)]}^{E^*}(g)(x)\big|+
\|P_{[1,n(m)]}^{E^*}(g)-g\|\big]\\
&\le (1+\ep)\big[\max_{g\in B_m}
|g\big(P_{[1,n(m)]}^{E}(x)\big)|+\ep\big]\le
(1+\ep)\big[\|P_{[1,n(m)]}^{E}(x)\|+\ep\big]\ ,
\end{align*}
which proves our claim and finishes the proof of part (a).

In order to show (b) we assume that  $X$ is a subspace of a
space $Z$ which has a boundedly complete FDD $(E_i)$ and the unit
ball of $X$ is a $w^*$-closed subset of $Z$, which is the dual of $\Zs$.

For $m\in\nat$ and $\ep>0$ we will show that the inequality in \eqref{E:2.4.1}
holds for some $n$ and $C=K$. 
Assuming that this was not true
we could choose a sequence $(x_n)\subset S_X$ so that for any $n\in\nat$
$$\|x_n\|_{X/X_m}>K\big[\|P^{E}_{[1,n]}(x_n)\|+\ep\big]\ .$$
By the compactness of $B_X$ in the $w^*$ topology we can choose
a subsequence $x_{n_k}$ which converges  $w^*$ to some $x\in B_X$.
For fixed $\ell$ it follows that $(P^{E}_{[1,\ell]}(x_{n_k}))$
converges in norm to $P^{E}_{[1,\ell]}(x)$. 
Secondly, since $X/X_m$ is finite dimensional it follows that
$\lim_{k\to\infty} \|x_{n_k}\|_{X/X_m}= \|x\|_{X/X_m}$,
and, thus, it follows that
\begin{align*}
 \|x\|&=\lim_{\ell\to\infty}\|P^{E}_{[1,\ell]}(x)\|\\
      &=\lim_{\ell\to\infty}\lim_{k\to\infty}
      \|P^{E}_{[1,\ell]}(x_{n_k})\|\\
       &\le K \limsup_{k\to\infty}\|P^{E}_{[1,n_k]}(x_{n_k})\| \\
        &\le \limsup_{k\to\infty} \|x_{n_k}\|_{X/X_m}-K\ep=
          \|x\|_{X/X_m}-K\ep\ ,
 \end{align*}
which is a contradiction since $\|x\|\ge \|x\|_{X/X_m}$.
\end{proof}

By combining Theorem~\ref{T:2.2} and Proposition~\ref{P:2.4} we deduce

\begin{cor}\label{C:2.6}
Let $X$ be a subspace of a space $Z$ with an FDD $(E_i)$ and assume that
this embedding satisfies condition \eqref{E:2.4.1}. 
Let $K\ge1$ be the projection constant
of $(E_i)$ in $Z$ and let $C\ge1$ be chosen so that   \eqref{E:2.4.1} holds.

For $\A\subset S^\omega_X$ the following  conditions are equivalent
\begin{enumerate}
\item[a)] For all null   sequences $\epb=(\ep_n)\subset (0,1]$,
 $WI(\,\overline{\A^X_{\epb}},Z)$ holds.
\item[b)]  For all null sequences $\epb=(\ep_n)\subset (0,1]$
there exists a blocking $(G_n)$ of $(F_n)$ so that 
every $\epb/420CK$-skipped block sequence $(z_n)\subset X$ 
of $(G_n)$ is in $\overline{\A_{\epb}}$.
\end{enumerate}
In the case that $X$ has a separable dual  (a) and (b) are
equivalent to the following condition
\begin{enumerate}
\item[c)] For all null sequences $\epb=(\ep_n)\subset (0,1]$ every weakly 
null tree in $S_X$ has a branch in $\overline{\A_{\epb}}$.
\end{enumerate}
In the case that  $(E_i)$ is a boundedly complete FDD of $Z$
and $B_{X}$ is $w^*$-closed in $Z=(Z^{(*)})^*$ the conditions (a) 
and (b) are equivalent to
\begin{enumerate}
\item[d)] For all null sequences $\epb=(\ep_n)\subset (0,1]$ every 
$w^*$-null tree in $S_X$ has a branch in $\overline{\A_{\epb}}$
\end{enumerate}
\end{cor}

\begin{proof}
(a) $\Rightarrow$ (b) Let $\epb=(\ep_i)\subset(0,1]$ be a null
sequence, choose $\etab=(\eta_i)$  with $\eta_i=\ep_i/3$, for $i\in\nat$,
and $\deltab=(\delta_i)$ with $\delta_i=\eta_i/140CK=\ep_i/420CK$.

We deduce from Proposition~\ref{P:2.4} that
$WI(\,\overline{\A_{\etab}}, (S^\omega_X)_{5\deltab})$ holds. 
Using Theorem~\ref{T:2.2} we
can block $(E_i)$ into $(G_i)$ so that every skipped block of $(G_i)$ in
$(S^\omega_X)_{5\deltab}$ (as a subset of $S_Z$) is in
$\overline{\A_{2\etab}}$ (actually we are using the quantified result 
given by the proof of Theorem~\ref{T:2.2}).

Assume $(x_i)\subset S_X$ is a $\deltab$-skipped block sequence of
$(G_i)$ and let  $1\le k_1< \ell_1< k_2< \ell_2<\cdots$ in $\nat$ so that
$$
\| x_n - P^E_{(k_{n},\ell_n]}(x_n)\|<\delta_n, \text{ for all }n\in\nat\ .$$
The sequence $(z_n)$ with
$z_n=P^E_{(k_{n},\ell_n]}(x_n)/\|P^E_{(k_{n},\ell_n]}(x_n)\|$, for $n\in\nat$,
is a skipped block sequence of $S_Z$ and we deduce that
\begin{align*}
\|x_n-z_n\|&\le\|x_n - P^E_{(k_{n},\ell_n]}(x_n)\|+
\|P^E_{(k_{n},\ell_n]}(x_n)\|\Big|1-\frac1{\|P^E_{(k_{n},\ell_n]}(x_n)\|}\Big|\\
&\le \delta_n+(1+\delta_n)\frac{\delta_n}{1-\delta_n}\le 5\delta_n\ .
\end{align*}
This implies that $(z_n)\in\overline{\A_{2\etab}}$ and thus 
by Lemma~\ref{L:2.1a} and our choice of $\etab$,
$$(x_i)\in \overline{(\,\overline{\A_{2\etab}})_{\etab}}\subset 
\overline{\A_{\epb}}\ ,$$
which finishes the verification of (b).

\noindent (b) $\Rightarrow$ (a) is clear since for any blocking
$(G_i)$ of $(E_i)$ and
 any null sequence $\deltab=(\delta_i)\subset(0,1]$ every tree
$(x_\alpha)_{\alpha\in T_\infty}$ in $S_X$ with the property that
$x_{(\alpha,n)}\in X_n$, whenever $n>n_\ell $
 and $\alpha=(n_1,\ldots, n_\ell)\in T_\infty\cup{\emptyset}$
 has a full subtree all of whose branches are
 $\deltab$-skipped block sequences of $(G_i)$.

Now assume that $X$ has a separable dual, or
$(E_i)$ is a boundedly complete FDD of $Z$ and  $B_X$ in $Z$ $w^*$-closed.

It is clear that (c) or (d), respectively, imply (a).
Secondly, since for any null sequence $\deltab=(\delta_i)\subset(0,1]$ and
any blocking $(G_i)$ every weakly null tree in $S_X$ 
(in the case that $X$, has  a separable dual)
or every $w^*$ null tree   (in the boundedly complete case) has a full
subtree  all of whose branches are 
$\deltab$-skipped block sequences of $(G_i)$ we deduce
that (b) implies (c) or (d) respectively.
\end{proof}

Motivated by the asymptotic structure of a Banach space 
we  introduce the following ``coordinate-free" variant of our games.
Again let $X$ be a separable Banach space and for $\A\subset
S_X^\omega$ we consider the following {\em coordinate-free} $\A$-{\em game}.
\begin{align*}
&\text{Player I chooses }X_1\in\cof(X)\ ,\\
&\text{Player II chooses } x_1\in X_1,\quad \|x_1\|=1\ ,\\
&\text{Player I chooses }X_2\in\cof(X)\ ,\\
&\text{Player II chooses } x_2\in X_2,\quad \|x_2\|=1\ ,\\
&\vdots
\end{align*}
As before, Player I  wins  if $(x_i)\in \A$. 
We will show that
$X$ can be embedded into a space $Z$ with an FDD so that
for all $\epb=(\ep_i)\subset(0,1]$ Player I has a winning
strategy in the coordinate-free 
$\overline{\A_{\epb}}$-game, which we will denote by
$WI(\,\overline{\A_{\epb}},\cof(X))$, 
if and only if for all  $\epb \subset(0,1]$
 he has a winning strategy for the $(\,\overline{\A_{\epb}},Z)$-game.

First note that since we only considering fattened sets  and their
closures, Player II has a winning strategy if and only if he has a
winning strategy choosing his vectors out of a dense and countable
subset of $S_X$ determined before the game starts. 
But this implies that there is countable set of
cofinite dimensional subspaces, say $\{Y_n:n\in\nat\}$
from which player I can choose if he has a winning strategy.
Moreover if we consider a countable set $\B$ of coordinate free
games, there is a countable set $\{Y_n:n\in\nat\}$ so that
for all $\A\in\B$
\begin{equation}\label{E:2.5} 
\forall\ \epb\subset(0,1]\quad  WI(\,\overline{\A_{\epb}},\cof(X))
   \iff  \forall\ \epb\subset(0,1)\quad
   WI(\,\overline{\A_{\epb}},\{Y_n:n\in\nat\})\ ,
 \end{equation}
where we write $ WI(\,\overline{\A_{\epb}},\{Y_n:n\in\nat\})$, if
 player I has a winning strategy for the coordinate-free $\A$-game,
even if he can only choose his spaces out of the set $\{Y_n:n\in\nat\}$.
Note that by passing to $(\bigcap_{i=1}^n Y_i)$ we can always
assume that the $Y_n$'s are decreasing in $n\in\nat$.
In case that $X$ has a separable dual and we let $(x^*_n)$  be a dense
 subset of $X^*$, we can put for $n\in\nat$
$$Y_n=\{x^*_1,x^*_2,\ldots, x^*_n\}^\perp 
=\{x\in X: \forall\  i\!\le\! n\,\, x^*_i(x)=0\}\ ,$$
and observe that \eqref{E:2.5} holds for all $\A\subset S^\omega_Z$.

The following result was shown in \cite[Lemma 3.1]{OS1} and its proof 
was based on techniques and results
of  W.B.~Johnson, H.~Rosenthal and M.~Zippin \cite{JRZ}. 

\begin{lem}\label{L:2.7}
Let $(Y_n)$ be a decreasing sequence of  closed subspaces of $X$,
each having finite codimension. 
Then $X$ is isometrically embeddable
into a space $Z$ having an FDD $(E_i)$ so that (we identify $X$ with
its isometric image in $Z$)
\begin{enumerate}
\item[a)]   $\coo(\oplus_{i=1}^\infty E_i)\cap X$ is dense in $X$.
\item[b)]  For every $n\in\nat$ the finite codimensional
subspace $X_n=\oplus_{i=n+1}^\infty E_i\cap X $ is contained in $Y_n$.
\item[c)]  There is a $c>0$, so that for every $n\in\nat$
there is a finite set $D_n\subset S_{\oplus_{i=1}^n E_i^*}$ such
that whenever $x\in X$
\begin{equation}\label{E4e}
\|x\|_{X/Y_n}=\inf_{y\in Y_n}\|x-y\| \le c\max_{w^*\in D_n} w^*(x) \ .
\end{equation}
\end{enumerate}
{From} (a) it follows that   $\coo (\oplus_{i=n+1}^\infty E_i)\cap
X$ is a dense linear subspace of $X_n$.

Moreover if $X$ has a separable dual $(E_i)$ can be chosen to be
shrinking (every normalized block sequence in $Z$ with respect to
$(E_i)$ converges weakly to $0$, or, equivalently, $Z^*
=\oplus_{i=1}^\infty E_i^*$), and if $X$ is reflexive $Z$ can also
be chosen to be reflexive.
\end{lem}

So assume that for a countable set $\B$ of games
that $(Y_n)$ is a sequence of decreasing finite codimensional
closed spaces satisfying the equivalences of \eqref{E:2.5}.
We then use Lemma~\ref{L:2.7} to embed
$X$ into a space $Z$ with an FDD  $(E_i)$.

Note that b) of Lemma~\ref{L:2.7} implies that for all $\A\in\B$
\begin{equation*}
\forall\ \epb\subset(0,1)\quad  WI(\,\overline{\A_{\epb}},\cof(X))
\iff  \forall\ \epb\subset(0,1)\qquad
WI(\,\overline{\A_{\epb}},Z)\ .
\end{equation*}
Using the embedding of $X$ given by Lemma~\ref{L:2.7}  a result
similar to Proposition~\ref{P:2.4} can be shown. 
The proof is very similar, therefore we will only present a sketch.

\begin{prop}\label{P:2.4a}
Assume that $X$ is a Banach space and $\{Y_n:n\in\nat\}$ a decreasing sequence
of cofinite dimensional subspaces. 
Let $Z$  be a space with an FDD
$(E_i)$ which satisfies the conclusion of Lemma~\ref{L:2.7}.

Assume that $\A\subset S_X^\omega$ such that
we have $WI(\,\overline{\A_{\epb}^X},\{Y_n:n\in\nat\})$  for all
null sequences  $\epb\subset(0,1]$.

Then  for all null sequences $\epb=(\ep_i)\subset (0,1]$,
$WI( \overline{\A_{\epb}},(S_X^\omega){\deltab})$ holds,
where $\deltab=(\delta_i)=(\ep_i/28cK)$, with  $c$ as in 
Lemma~\ref{L:2.7}, $K$ is the projection constant of $(E_i)$ in $Z$, and
where the fattenings $\overline{\A_{\epb}}$ and
$(S_X^\omega){\deltab}$ are taken in $Z$.
\end{prop}

\begin{proof}[Sketch of  proof] 
Note that instead of condition
\eqref{E:2.4.1} the following condition is satisfied.
\begin{align}\label{E:2.4a.1} 
&\text{There is a $C>0$ so that for all $m\in\nat$}\\
&\qquad\|x\|_{X/Y_m}\le C \|P^E_{[1,m]}(x)\| \text{ whenever }x\in
 S_X.\notag
\end{align}
Also note that $WI(\,\overline{\A_{\epb}},\{Y_n:n\in\nat\})$ means that
every tree $(x_\alpha)\subset S_X$, with  the property
that for $\alpha=(m_1,m_2,\ldots, m_\ell)\in T_\infty$
we have that $x_\alpha\in Y_{m_\ell}$, has a branch in $\overline{\A_{\epb}}$.

We follow the proof of Proposition~\ref{P:2.4} until choosing the
$x_\alpha$'s which we will not choose in $X_{m_\ell}$ but in
$Y_{m_\ell}$ instead. 
Then the proof continues as the proof of Proposition~\ref{P:2.4}.\end{proof}

Using Proposition~\ref{P:2.4a} and Theorem~\ref{T:2.2} we deduce
the following.

\begin{cor}\label{C:2.8}
Let $\B$ be a countable set of $\A\subset S_X^\omega$
and assume that $Z$ is a space with an $\FDD\, (E_i)$ which contains $X$ and
satisfies the conclusion of Lemma~\ref{L:2.7}.

For $\A\in \B$ the following  conditions are equivalent:
\begin{enumerate}
\item[a)] For all null   sequences $\epb=(\ep_n)\subset (0,1]$,
$WI(\,\overline{\A^X_{\epb}},\cof)$ holds.
\item[b)] For all null   sequences $\epb=(\ep_n)\subset (0,1]$,
$WI(\,\overline{\A^X_{\epb}},Z)$ holds.
\item[c)]  For all null sequences $\epb=(\ep_n)\subset (0,1]$,
there exists a blocking $(G_n)$ of $(E_n)$ so that 
every $\epb/420CK$-skipped block sequence $(z_n)\subseteq X$ 
of $(G_n)$ is in $\A_{\epb}$.
\end{enumerate}
In the case that $X$ has a separable dual  (a), (b) and (c) are
equivalent to the following condition (which is independent of the 
choice of $Z$).
\begin{enumerate}
\item[d)] For all null sequences $\epb=(\ep_n)\subset (0,1]$ every 
weakly null tree in $S_X$ has a branch in $\overline{\A_{\epb}}$.
\end{enumerate}
Moreover, in the case that $X$ has a separable dual
we deduce from the remarks after the equivalence \eqref{E:2.5}, 
Corollary~\ref{C:2.6} and Proposition~\ref{P:2.5} that above 
equivalences hold for {\em any embedding }of $X$
into a space $Z$ having a shrinking FDD.
\end{cor}

\section{Banach spaces of Szlenk index $\w$}        

In this section we will present (Theorem~\ref{thm:equivconditions})  
a long list of equivalent conditions for a space $X$ to have Szlenk index $\w$.
We also show how Kalton's $c_0$  theorem
(Theorem~\ref{thm:Kalton})      
can be proved with our techniques. 
We begin with some definitions
that will be used in later sections as well as this one.

\begin{Defns}
Let $1\le q\le p\le\infty$ and $C<\infty$. 
A (finite or infinite) $\FDD$ $(E_i)$ for a
Banach space $Z$ is said to satisfy {\em $C$-$(p,q)$-estimates\/} if
for all $n\in\nat$ and block sequences $(x_i)_{i=1}^n$ w.r.t.\ $(E_i)$,
$$C^{-1} \bigg( \sum_1^n \|x_i\|^p\bigg)^{1/p}
\le \Big\| \sum_1^n x_i\Big\| \le C\bigg( \sum_1^n
\|x_i\|^q\bigg)^{1/q}\ .$$

A space $X$ satisfies {\em $C$-$(p,q)$-tree estimates\/} if for all
weakly null trees in $S_X$ there exist branches $(x_i)_{i=1}^\infty$
and $(y_i)_{i=1}^\infty$ satisfying for all $(a_i) \in c_{00}$,
\begin{equation}\label{eq3.1}
C^{-1} \Big(\sum |a_i|^p\Big)^{1/p} \le \Big\| \sum a_i x_i\Big\| \
\text{ and }\ \Big\| \sum a_i y_i\Big\|  \le C\Big( \sum
|a_i|^q\Big)^{1/q}\ .
\end{equation}

If $X\subseteq Y^*$, a separable dual space, we say that $X$
satisfies {\em $C$-$(p,q)$-$w^*$-tree estimates\/} if each $w^*$
null tree in $S_X$ admits branches $(x_i)$ and $(y_i)$ satisfying
\eqref{eq3.1}.

We will say that $X$ satisfies {\em $(p,q)$-tree estimates\/} if it
satisfies $C$-$(p,q)$-tree estimates for some $C<\infty$ and
 similarly for $(p,q)$-$w^*$ tree estimates.
\end{Defns}

It is perhaps worth noting that if every weakly null tree in $X$
admits a branch dominating the unit vector basis of $\ell_p$
(not assuming that the constant of domination can be chosen independently
of the tree) then
$X$ satisfies $(p,1)$-tree estimates (and similar remarks hold for
$(\infty,q)$-tree estimates or $(p,q)$-$w^*$-tree estimates).
Indeed, if no uniform constant existed one could assemble a tree
with no branch dominating the unit vector basis of $\ell_p$
\cite{OSZ}[Proposition 1.2].

\begin{Defn} \cite{Pr}
Let $1\le p\le\infty$ and let $Z$ be a Banach space with an $\FDD E= (E_i)$. 
Then $Z_p (E)$ is the completion of $c_{00} (\oplus_{i=1}^\infty E_i)$ under
$$\|z\|_p = \sup \bigg\{ \Big( \sum_j \| P_{I_j}^E z\|^p\Big)^{1/p} :
I_1 < I_2 <\cdots \text{ are intervals in } \nat\bigg\}\ .$$
$(E_i)_{i=1}^\infty$ is then a bimonotone $\FDD$ for $Z_p(E)$ which
satisfies 1-$(p,1)$-estimates. 
Moreover, if $Z$ is isomorphic to
$\tZ$ then $Z_p(E)$ is naturally isomorphic to $\tZ_p(E)$.
\end{Defn}

Our main tool for proving Theorem~\ref{thm:equivconditions} is
 the following result which is a non-reflexive version of Theorem~2.1~a) in
\cite{OS2}.

\begin{thm}\label{thm:main}
Let $Z$ be a Banach space with a boundedly complete $\FDD$ $E=
(E_i)$ and let  $X$ be a subspace of $Z$ with $B_X$ being a
$w^*$-closed subset of $Z$ ($=(Z^{(*)})^*$). 
Let $1\le p\le \infty$. 
If $X$ satisfies $(p,1)$-$w^*$-tree estimates in $Z$ then
there exists a blocking $F= (F_i)$ of $(E_i)$ so that $X$ naturally
embeds into $Z_p(F)$.
\end{thm}

To prove this we need a blocking lemma which appears in various
forms in \cite{KOS}, \cite{OS1}, \cite{OS2}, \cite{OSZ} and
ultimately results from a blocking trick of Johnson \cite{J1}. 
We will use this lemma as well in section \ref{S:4}.

\begin{lem}\label{lem:blocking}
Let $X$ be a  subspace of a space $Z$ 
having a boundedly complete $\FDD$ $(E_i)$ with projection constant $K$
with $B_X$ being a $w^*$-closed subset of $Z$.     
Let $\delta_i\downarrow 0$.
Then there exists a blocking $(F_i)$ of $(E_i)$ given by
$F_i = \oplus_{j = N_{i-1}+1}^{N_i} E_j$ for some
$0= N_0 < N_1<\cdots$  with the following properties.
For all $x\in S_X$ there exists $(x_i)_{i=1}^\infty \subseteq X$
and for all $i\in \nat$ there exists $t_i\in (N_{i-1},N_i)$ satisfying
($t_0=1$ and $t_1>1$)
\begin{itemize}
\item[a)] $x = \sum_{j=1}^\infty x_j$,
\item[b)] $\|x_i\| < \delta_i$ or $\|P_{(t_{i-1},t_i)}^E x_i - x_i\|
<\delta_i \|x_i\|$,
\item[c)] $\|P_{(t_{i-1},t_i)}^E x-x_i\| < \delta_i$,
\item[d)] $\|x_i\| < K+1$,
\item[e)] $\|P_{t_i}^E x\| < \delta_i$.
\end{itemize}
Moreover, the above hold for any blocking of $(F_i)$
(which would redefine the $N_i$'s).
\end{lem}

\begin{proof}
We observe that for all $\ep>0$ and $N\in\nat$ there exists $n>N$
 such that if $x\in B_X$, $x= \sum y_i$ with $y_i\in E_i$ for all $i$,
 then there exists $t\in (N,n)$ with
$$\|y_t\| <\ep\ \text{ and }\ \dist \bigg( \sum_{i=1}^{t-1} y_i,X\bigg)
 < \ep\ .$$
Indeed, if this was not true  for any $n>N$ we can find 
$y^{(n)}\in B_X$ failing the conclusion for $t\in (N,n)$.
Choose a subsequence of $(y^{(n)})$ converging  $w^*$  to $y\in X$
and choose $t>N$ so that $\|P_{[t,\infty)}^E y\| < \ep/2K$.
Then choose $y^{(n)}$ from the subsequence so that $t<n$ and
$\| P_{[1,t]}^E (y-y^{(n)})\| < \ep/2K$.
Thus
$$\| P_{[1,t)}^E y^{(n)} - y\|
\le \| P_{[1,t)}^E (y^{(n)} - y)\| + \| P_{[t,\infty)}^E y\|
 < \frac{\ep}{2K} + \frac{\ep}{2K} <\ep\ .$$
Also
$$\| P_t^E y^{(n)}\|
\le\| P_t^E (y^{(n)}-y)\|  + \| P_t^E y\| < \frac{\ep}2 +\frac{\ep}2 =\ep\ .$$
This contradicts our choice of $y^{(n)}$.

Let $\ep_i\downarrow0$ and by the observation choose $0= N_0 < N_1
<\cdots$ so that for all $x\in S_X$ there exists $t_i\in (N_{i-1},N_i)$
and $z_i \in X$ with $\|P_{t_i}^E x\| <\ep_i$ and $\|P_{[1,t_i-1)}^E
x-z_i\| < \ep_i$ for all $i\in \nat$.
Set $x_1 = z_1$ and $x_i = z_i - z_{i-1}$ for $i>1$.
Thus $\sum_{i=1}^n x_i = z_n \to x$ so a) holds.
Also
$$\|P_{(t_{i-1},t_i)}^E x-x_i\|
 \le \|P_{[1,t_i)}^E x-z_i\|
 + \| P_{[1,t_{i-1}]}^E x-z_{i-1}\|
 < \ep_i + 2\ep_{i-1}\ ,$$
and
$$\| P_{(t_{i-1},t_i)}^E x_i - x_i\|
= \| (I- P_{(t_{i-1},t_i)}^E ) (x_i - P_{(t_{i-1},t_i)}^E x)\|
< (K+1) (\ep_i + 2\ep_{i-1})\ .$$
{From} these inequalities b), c) and d) follow if we take
$(\ep_i)$ so that $(K+1)(\ep_i + 2\ep_{i-1}) <\delta_i^2$.
\end{proof}

\begin{proof}[Proof of Theorem~\ref{thm:main}]
We may assume that $E$ is a bimonotone $\FDD$ for $Z$ and
that $X$ satisfies $C$-$(p,1)$-$w^*$-tree estimates in $Z$.

Let
\begin{equation*}
\A =\Big\{(x_i)_{i=1}^\infty \in S_X^\w :
\Big\|\sum a_i x_i\Big\| \ge C^{-1} \Big(\sum |a_i|^p\Big)^{1/p}
\text{ for all }(a_i)\in c_{00}\Big\}\ .
\end{equation*}
Choose a null sequence $\epb=(\ep_i)\subset(0,1)$ so that
\begin{equation*}
\overline{A_{\epb}}   \subset \Big\{ (x_i)_{i=1}^\infty \in S_X^\w :
\Big\|\sum a_i x_i\Big\| \ge (2C)^{-1} \Big(\sum |a_i|^p\Big)^{1/p}
\text{ for all } (a_i)\in c_{00}\Big\}\ .
\end{equation*}
By Corollary~\ref{C:2.6}
there exist $\bar\delta  = (\delta_i)$ with $\delta_i
\downarrow 0$ and a blocking of $(E_i)$, which we still denote by
$(E_i)$ so that if $(x_i)_{i=1}^\infty \in S_X^\w$ is a
$\bar\delta$-skipped block sequence w.r.t.\ $(E_i)$ then $(x_i)\in
\bA_{\epb}$. Wlog $\sum_1^\infty \delta_i < \frac12$. 
We will produce a blocking $(F_i)$ of $(E_i)$ and $A<\infty$ so that for all
$0=n_0 < n_1 <\cdots$ and  $x\in S_X$,
\begin{equation}\label{eq3.2}
\bigg( \sum_{j=1}^\infty \| P_{(n_{j-1},n_j]}^F x\|^p \bigg)^{1/p}
\le A
\end{equation}
and this will finish the proof.

$(F_i)$ will be the blocking given by Lemma~\ref{lem:blocking} 
for $(\delta_i)$. 
We will show, using Lemma~\ref{lem:blocking}, that
\begin{equation}\label{eq3.3}
\bigg(\sum_j \|P_j^F x\|^p\bigg)^{1/p} \le A
\end{equation}
and by the ``moreover'' part of the lemma the same proof will yield
\eqref{eq3.2}.

Let $F_j = \oplus_{i=N_{j-1} +1}^{N_j} E_i$ be as in  
Lemma~\ref{lem:blocking}, 
$x\in S_X$ and let $(x_i)$, $(t_i)$ be as in the Lemma. 
Set $B = \{ i \ge 2$; $x_i\ne0$ and 
$\|P_{(t_{i-1},t_i)}^E x_i - x_i\| < \delta_i \|x_i\|\}$ 
and for $i\in B$ let $\bar x_i = x_i/\|x_i\|$. 
Note that if $i\ge 2$, $i\notin B$, then $\|x_i\| < \delta_i$.

Now $(\bar x_i)_{i\in B}$ is a $\bar \delta$-skipped block sequence
w.r.t.\ $(E_i)$ and so
$$2C\Big\| \sum_{i\in B} x_i \Big\| \ge \bigg( \sum_{i\in B}
\|x_i\|^p \bigg)^{1/p}\ .$$ 
Also
$$\Big\| \sum_{i\notin B} x_i\Big\|
\le \|x_1\| + \sum \delta_i < 2+1 = 3\qquad (\text{since } K=1)\ .$$
Thus
$$\Big\| \sum_{i\in B} x_i\Big\|  \le \|x\| + \Big\| \sum_{i\notin B}
x_i\Big\| <4\ .$$

It follows that
\begin{gather*}
\sum \| x_i\|^p \le \|x_1\|^p + \sum_{i\in B} \|x_i\|^p
+ \sum_{\substack{i\notin B\\ i\ge 2}} \|x_i\|^p\\
< 2^p + (8C)^p + 1 \equiv D^p\ .
\end{gather*}
For $i\in\nat$ set $y_i = P_{(t_{i-1},t_i]}^E x$. 
Then
$$\| y_i - x_i\|
\le \| P_{(t_{i-1},t_i)}^E x-x_i\| + \| P_{t_i}^E x\| < 2\delta_i\ .$$ 
Hence
$$\Big(\sum \|y_i\|^p\Big)^{1/p}
\le \Big(\sum \|x_i\|^p\Big)^{1/p} + \Big(\sum
(2\delta_i)^p\Big)^{1/p} < D+1\ .$$ 
Finally write $x= \sum z_i$
where $z_i \in F_i$ for all $i$. 
Then $z_i = P_i^F (y_i + y_{i+1})$, so $\|z_i\| \le \| y_i\| + \| y_{i+1}\|$. 
Hence
$$\Big(\sum \|z_i\|^p\Big)^{1/p} \le 2 (D+1) \equiv A\ .$$
\end{proof}

Before stating Theorem~\ref{thm:equivconditions} we need some
definitions and preliminaries.

\begin{Defn}
$X$ has the {\em $w^*$-UKK} if for all $\ep>0$ there exists
$\delta^* (\ep) >0$ so that if $(x_n^*) \subseteq B_{X^*}$ converges
$w^*$ to $x^*$ and $\liminf_{n\to\infty} \|x_n^* - x^*\| \ge \ep$
then $\|x^*\| \le 1- \delta^* (\ep)$.
\end{Defn}

We have defined $\Sz(X)$ but  there is another view of this
index which shall prove useful.
In \cite[Theorem~4.2]{AJO} an index $I_{\ell_1^+,w}(X) $ is defined and  shown
to equal $\Sz(X)$ if $X$ does not contain an isomorph of $\ell_1$.
The precise definition need not concern us here.
However we note that one consequence is $\Sz(X) =\w$
iff $\forall\ K >1\ \exists\ n(K)$ so that if $(e_i)_{i=1}^n \in \{X\}_n$
is an $\ell_1^+ - K$ sequence, i.e., $\|\sum_1^n a_i e_i\| \ge
K^{-1} \sum_{i=1}^n a_i$ for $(a_i)_{i=1}^n \subseteq [0,\infty)$,
then $n\le n(K)$.

\begin{Defn}
For a Banach space $X$ the {\em modulus of  asymptotic
uniform smoothness\/} $\bar\rho_X(t)$ is given for $t>0$ by
$$\bar \rho_X(t) = \sup_{\|x\| =1} \inf_{Y\in\cof(X)} \sup_{y\in tB_Y}
 \|x+y\| -1\ .$$

The modulus of asymptotic uniformly convexity $\bar\delta_X(t)$ is
given for $t>0$ by
$$\bar\delta_X(t) = \inf_{\|x\|=1} \sup_{Y\in\cof(X)}
 \inf_{\substack{y\in Y\\ \|y\| \ge t}} \|x+y\| -1\ .$$

$X$ is {\em asymptotically uniformly smooth} (a.u.s.) if
$\lim_{t\to 0^+} \bar \rho_X(t)/t = 0$.

$X$ is {\em asymptotically uniformly convex\/} (a.u.c.) if for
$t>0$, $\bar \delta_X (t) >0$.

$X$ is {\em a.u.s.\ of power type $p$} if for some $K<\infty$,
$\bar \rho_X (t) \le K\,t^p$ for $t>0$.

$X$ is {\em a.u.c.\ of power type $p$} if for some $K>0$,
$\bar \delta_X (t) \ge K\,t^p$ for $t>0$.

If $X$ is a dual space we can define similar modulii $\bar\delta_X^*(t)$
and $\bar\rho_X^* (t)$ using
\begin{equation*}
\cof^*(X)  = \Big\{ Y\subseteq X : Y\text{ is a $w^*$-closed finite}
\text{ co-dimensional subspace of } X\Big\}\ .
\end{equation*}
\end{Defn}

More about these modulii can be found in \cite{JLPS} but we shall
extract a few things we need in proving Theorem~\ref{thm:equivconditions}.
A.u.s.\ and a.u.c.\ say something
about weakly null trees and $\{X\}_n$.
Let $(e_i)_{i=1}^n \in \{X\}_n$ and let $(a_i)_{i=1}^n \subseteq (0,1]$.
Assume that $\bar\delta_X(t) \ge K\,t^p$ for some $K$ and all $t>0$.
Using that there exists $c>0$ with $K\,t^p \ge (1+c\,t^p)^{1/p}-1$
for $t>0$ we obtain
\begin{equation*}
\begin{split}
\Big\|\sum_1^n a_i e_i\Big\|^p
& \ge \Big\| \sum_1^{n-1} a_i e_i\Big\|^p
\bigg(1+\bar\delta_X 
\Big(\frac{|a_n|}{\|\sum_1^{n-1} a_i e_i\|}\Big)\bigg)^p\\
\noalign{\vskip6pt}
&\ge \Big\| \sum_1^{n-1} a_i e_i\Big\|^p + c|a_n|^p\\
\noalign{\vskip6pt}
&\ge \cdots \ge c\sum_1^n |a_i|^p\ .
\end{split}
\end{equation*}
Similarly if we begin with a weakly null tree in $S_X$ we can extract
a branch $(x_i)$ satisfying
$$\Big\| \sum a_i x_i\Big\| \ge \frac{c}2 \Big( \sum |a_i|^p\Big)^{1/p}$$
for all $(a_i)\in c_{00}$.

With a similar argument for $\bar\rho_X (t)$ and $\bar\sigma_X^*(t)$ we obtain

\begin{prop}\label{prop:JLPS}\cite{JLPS}
\begin{itemize}
\item[a)]  
If $X$ is a.u.c.\ of power type $p$ then $X$ satisfies $(p,1)$-tree
estimates.
\item[b)]  
If $X$ is a.u.s.\ of power type $q$ then $X$ satisfies
$(\infty,q)$-tree estimates.
\item[c)] Let $X= Y^*$ be a separable dual.
If $X$ is $w^*$-a.u.c.\ of power type $p$ (i.e., $\bar\sigma_X^* (t)
\ge K\,t^p$) then $X$ satisfies $(p,1)$-$w^*$-tree estimates.
\end{itemize}
\end{prop}

\begin{thm}\label{thm:equivconditions}   
Let $X^*$ be separable.
The following are equivalent.
\begin{itemize}
\item[(1)] $\Sz(X) =\w$.
\item[(2)] $\exists\ q > 1\ \exists\  K<\infty\ \forall\ n\ \forall\
(e_i)_{i=1}^n \in \{X\}_n\ \forall\ (a_i)_{i=1}^n \subseteq \real$,
$$\Big\| \sum_1^n a_i e_i\Big\| \le K\Big( \sum_1^n |a_i|^q\Big)^{1/q}\ .$$
\item[(3)] $\exists\ q>1$ so that $X$ satisfies $(\infty,q)$-tree estimates.
\item[(4)] $\exists\ p<\infty$ so that $X^*$ satisfies $(p,1)$-$w^*$-tree
estimates.
\item[(5)] $\exists\ p<\infty\ \exists$ a Banach space $Z$ with a
boundedly complete $\FDD\ E$ so that $X^*$ embeds into $Z_p(E)$ as a
$w^*$-closed subspace.
\item[(6)] $X$ can be renormed to be a.u.s.\ of power type $q$
for some $q>1$.
\item[(7)] $X$ can be renormed to be a.u.s.
\item[(8)] $X$ can be renormed so that $\bar\rho_X(t) <t$ for some
$t>0$.
\item[(9)] $X$ can be renormed to be $w^*$-UKK with modulus
$\delta^* (\ep) \ge c\,\ep^p$ for some $p<\infty$.
\item[(10)] $X$ can be renormed to be $w^*$-UKK.
\item[(11)] $\exists\ p<\infty$ so that $X$ can be renormed so that
$\bar \delta_X^* (t)$ is of power type $p$.
\item[(12)] $X$ can be renormed to be $w^*$-a.u.c.
\end{itemize}
\end{thm}

\begin{proof}

\noindent $(2) \Rightarrow (1)$. (2) implies that $I_{\ell_1^+,\w}
(X) = \w = S_z(X)$ by our earlier remarks.
\smallskip

\noindent $(1)\Rightarrow (2)$. This follows from the fact that for
$n\in\nat$ there exists $q>1$ so that every normalized monotone
basis which does not admit a normalized block basis of length $n$
which is $\ell_1^+$ with constant~2 is $6$-dominated by the unit
vector basis of $\ell_q^n$ (proved in \cite{Ja}, \cite{J2}).
Since $I_{\ell_1^+,w}(X)=\omega$ (2) follows by our earlier remarks
and the fact that if $(x_i)_{i=1}^m$ is a normalized block 
sequence of some sequence
$(e_i)_{i=1}^n\in\{X\}_n$ then $(x_i)_{i=1}^m\in \{X\}_m$.

\smallskip

\noindent $(2)\Rightarrow (3)$. Let $X\subseteq Y$, a space with a
shrinking $\FDD$ $(E_i)$ (by 1.2). Using our discussion of
asymptotic structure, applying  Corollary~\ref{C:2.8} to
$\B=\{\A^{(n)}:n\in\nat\}$, with
\begin{equation}
\A^{(n)}=\left\{ (x_i)\in S_X^\omega:
\begin{matrix}\exists\  (e_i)_{i=1}^n \in\{X\}_n\text{ so that}\\
\noalign{\vskip6pt}
(x_i)_{i=1}^n\text{ is $2$-equivalent to}(e_i)_{i=1}^n
\end{matrix}\right\}
\end{equation}
 and a diagonal argument we can find
$\bar\delta = (\delta_i)$, $\delta_i\downarrow 0$, and a blocking
$(F_i)$ of $(E_i)$ with the following property. For all $n\in\nat$
if $(x_i)_{i=1}^n \subseteq S_X$ is a
$(\delta_i)_{n+1}^{2n}$-skipped block sequence w.r.t.\
$(\oplus_{i=1}^n F_i, F_{n+1},F_{n+2},\ldots)$ then
$d_b((x_i)_{i=1}^n,\{X\}_n)<2$. Let $(x_\alpha)_{\alpha\in
T_\infty}$ to be a weakly null tree in $X$. 
Then the exists a branch
$(x_i)_{i=1}^\infty$ so that for all $n$ if $(y_i)_{i=1}^n$ is a
normalized block basis of $(x_i)_n^\infty$ then $(y_i)_{i=1}^n$ is a
$(\delta_i)_{n+1}^{2n}$-skipped block sequence w.r.t.\ $(\oplus_1^n
F_i, F_{n+1},\ldots)$ and so satisfies $2K$-upper $\ell_q^n$ estimates. 
Now it follows \cite[Proposition 3.5]{KOS} that for any
$q>\bar q>1$, $(x_i)$ satisfies $(\infty,\bar q)$-estimates. 
Thus (3) holds.
\renewcommand{\qed}{}
\end{proof}

\noindent $(3)\Rightarrow (4)$ follows from the following

\begin{lem}\label{lem3.5}
Let $X^*$ be separable. 
If $X$ satisfies $(\infty,q)$-tree estimates for $q>1$ then $X^*$ 
satisfies $(q',1)$-$w^*$-tree estimates $(\frac1q + \frac1{q'} =1)$.
\end{lem}

\begin{proof}
Let $X$ satisfy $K$-$(\infty,q)$-tree estimates. 
Note the following.
If $(x_i^*)$ is normalized $w^*$-null in $X^*$ then there exists
$(x_i)\subseteq S_X$, $(x_i)$ is weakly null, and a subsequence
$(x_{n_i}^*)$ of $(x_i^*)$ with $\lim_i x_{n_i}^* (x_i) \ge \frac12$. 
Indeed we choose $(y_i)\subseteq S_X$ with $\lim x_i^*
(y_i)=1$ and pass to a weak Cauchy subsequence $(y_{k_i})$ so that
$\lim_i x_{k_i}^* (y_{k_{i-1}})=0$. 
Let $x_{n_i}^* = x_{k_{2i}}^*$
and $x_i = (y_{k_{2i}} - y_{k_{2i-1}}) /\|y_{k_{2i}}- y_{k_{2i-1}}\|$.

Let $(x_\alpha^*)_{\alpha\in T_\infty}$ be a $w^*$-null tree in $X^*$. 
Using the above remark we can pass to a full subtree 
which we still denote by $(x_\alpha^*)_{\alpha\in T_\infty}$ and find
a weakly null tree $(x_\alpha)_{\alpha\in T_\infty}\subseteq S_X$ so
that $x_\alpha^* (x_\alpha) > 1/3$ for all $\alpha$. 
By further pruning we can also assume that, given $\eta >0$,
$|x_\alpha^*(x_\beta)| <2^{-m-n}\eta$ and $|x_\beta^* (x_\alpha)| <
2^{-m-n}\eta$ if $\alpha <\beta$ and $|\alpha| = m$, $|\beta |=n$.
This pruning uses only that each node in $(x_\alpha^*)$ is
$w^*$-null and each node in $(x_\alpha)$ is weakly null. 
An easy calculation shows that if $(x_i)_{i=1}^\infty$ is a branch in
$(x_\alpha)$ which is $K$-dominated by the unit vector basis of
$\ell_q$, then the corresponding  branch $(x_i^*)$ in $(x_\alpha^*)$
satisfies, for small $\eta$,
$$\Big\| \sum a_i x_i^*\Big\| \ge \frac1{3K} - \frac1K \eta > \frac1{4K}$$
if $(\sum |a_i|^{q'})^{1/q'} =1$.\qed
\smallskip

\noindent $(4) \Rightarrow (5)$. By  1.1 $X$ is a quotient of a space
with a shrinking basis and hence $X^*$ embeds as a $w^*$-closed
subspace into a space $Z$ with a boundedly complete FDD $E= (E_i)$.
Since any $w^*$-null tree in $S_{X^*}$ is a $w^*$-null tree w.r.t.\
$Z$, (5) follows from (4) by Theorem~\ref{thm:main}.
\smallskip

\noindent $(5)\Rightarrow (6)$. 
Let $X^*$ be embedded into $Z_p(E)$ as in (5). 
We renorm $X$ via
$$|x| = \sup \{ |x^* (x)| : x^* \in X^*,\ \|x^*\|_p\le 1\}\ .$$
It follows easily that $\bar\rho_X(t) \le (1+t^q)^{1/q} -1$ where
$\frac1p +\frac1q =1$ which proves (6).
\smallskip

\noindent $(6)\Rightarrow (7)\Rightarrow (8)$ is trivial.
\smallskip

\noindent $(8)\Rightarrow (1)$. 
Assume (1) fails so $\Sz(X) = I_{\ell_1^+,w} (X) >\w$. 
Then there exists $K\ge1$ so that for all
$n$ there exists an $\ell_1^+ -K$ sequence in $\{X\}_n$. 
By James' argument that $\ell_1$ is not distortable (which also works in the
$\ell_1^+$ case) we obtain that there exists $(e_1,e_2) \in \{X\}_2$
with
$$\|e_1 + te_2\| = 1+t \ \text{ for all }\ t>0\ .$$
Since $\Sz(X)$ is an isomorphic invariant, we have for all
renormings of $X$, $\bar\rho_X(t)=t$ for all $t>0$. 
Thus (8) fails.
\smallskip

\noindent $(5)\Rightarrow (9)$ by the renorming used in
$(5)\Rightarrow (6)$. 
Indeed if $(x_n^*) \subseteq
S_{X^*,\|\cdot\|_p}$ with $x_n^* \xrightarrow{w^*} x^*$ and $\lim_n
\|x_n^* - x^*\|_p \ge \ep$ then $\|x^*\|^p + \ep^p \le 1$.
\smallskip

\noindent $(9)\Rightarrow (10)\Rightarrow (1)$ is trivial.
\smallskip

\noindent $(5)\Rightarrow (11)$ holds again by the $(5)\Rightarrow
(6)$ argument.

\noindent $(11)\Rightarrow (12)$ is trivial.

\noindent $(12)\Rightarrow (4)$ 
Assume (12) holds. 
By \cite{GKL} (7) holds. 
Alternatively, it follows that there exists $n_0\in\nat$ so
that if $(e_i)_{i=1}^n$ is in the $w^*$-asymptotic structure of the
$w^*$-a.u.c.\ space $X^*$ and there exist $(a_i)_{i=1}^n\subseteq
[\frac12,1]$ with $\|\sum_1^n a_i e_i\|\le 1$ then $n\le n_0$.
Indeed we obtain $\| \sum_1^n a_i e_i\| \ge \frac12 [1+\deltab_{X}^*
(\frac12)]^{n-1}$. 
This condition yields that there exists $p = p(n_0) <\infty$ so that 
the unit vector basis of $\ell_p^n$
2-dominates $(e_i)_{i=1}^n$ for all $n\in\nat$, 
(\cite{Ja}, \cite{J2}, \cite{KOS}). 
Arguing then as in $(2)\Rightarrow (3)$ we obtain (4).
\end{proof}

We end this section with Kalton's $c_0$-theorem.

\begin{Defn}
$X$ has the {\em bounded tree property\/} if there exists $C<\infty$
so that for all weakly null trees in $S_X$ there exists a branch
$(x_i)_{i=1}^\infty$ with
$$\sup _n \Big\| \sum_{i=1}^n x_i\Big\| \le C\ .$$
Note that if $X$ has the bounded tree property and does not contain
an isomorph of $\ell_1$ then $\Sz(X) = I_{\ell_1^+,\w} (X) =\w\ .$
\end{Defn}

\begin{thm}\label{thm:Kalton} 
\cite{K}  Let $X$ have the bounded tree property. 
If $X$ does not contain an isomorph of $\ell_1$, then $X$ embeds into $c_0$.
\end{thm}

\begin{proof}
By (1.2)  we may regard $X\subseteq Z$, a space with a bimonotone
shrinking FDD $E= (E_i)$. 
Assume that $X$ has the bounded tree property with constant $C$. 
Let
$$\A = \bigg\{ (x_i)_{i=1}^\infty \in S_X^\w :
\sup_n \Big\| \sum_1^n x_i\Big\| \le C\bigg\}\ .$$ 
Choose $\epb \subseteq (0,1)$  so that 		
$$\bar\A_{\epb} \subseteq  \bigg\{ (x_i) \in S_X^\w :
\sup_n \Big\| \sum_1^n x_i\Big\| \le 2C\bigg\}\ .$$
By Corollary~\ref{C:2.6} we may choose $\bar \delta = (\delta_i)$,
$\delta_i\downarrow 0$, and a blocking of $E$ which we still denote
by $E= (E_i)$ so that any $\bar\delta$-skipped block sequence
$(x_i)\subseteq S_X$ w.r.t. $(E_i)$ is in $\bar\A_{\epb}$. 
Since $(\pm x_i)$ is a $\bar\delta$-skipped block sequence when $(x_i)$ is a
$\bar\delta$-skipped block sequence it follows by a convexity
argument that $\|\sum a_i x_i\| \le 2C$ for $(a_i)\in c_{00}$,
$(a_i)\subseteq [-1,1]$.

It follows that $X$ satisfies $(\infty,\infty)$-tree estimates and
hence by Lemma~\ref{lem3.5}, $X^*$ satisfies $(1,1)$-$w^*$-tree estimates. 
By Theorem~\ref{thm:main}, $X^*$ embeds as a $w^*$ closed
subspace into some space $Z^*_1(F_i^*)$ which is
$(\oplus_{i=1}^\infty F_i^*)_{\ell_1}$,
 where $F^* = (F_i^*)$ is some blocking of $(E_i^*)$.
{From} basic functional analysis we have that $X$ is a quotient of
$(\sum F_i)_{c_0}$. Hence $X$ is isomorphically a subspace of a
quotient of $c_0$ and hence embeds into $c_0$ since every quotient
of $c_0$ embeds into $c_0$. \cite{JZ1}.
\end{proof}

\section{Reflexive spaces}\label{S:4}

In this section we  first discuss the problem of characterizing when
a reflexive space $X$ satisfies $(p,q)$-tree estimates for a given
$1\le q\le p\le\infty$. 
The ultimate result is

\begin{thm}\label{thm4.1}
Let $X$ be a reflexive Banach space and let $1\le q\le p\le\infty$.
The following are equivalent
\begin{itemize}
\item[a)] $X$ satisfies $(p,q)$-tree estimates.
\item[b)] $X$ embeds into a reflexive space $Z$ having an FDD which satisfies
$(p,q)$-estimates.
\item[c)] $X$ is isomorphic to a quotient of a reflexive space $Z$ having an FDD
which satisfies $(p,q)$-estimates.
\item[d)] $X^*$ satisfies $(q',p')$-tree estimates where $1/q' + 1/q=1$ and
$1/p' + 1/p=1$.
\item[e)] $X$ embeds into a reflexive space $Z$ having an FDD which satisfies
1-$(p,q)$-estimates.
\end{itemize}
\end{thm}

The duality between an FDD $(E_i)$ satisfying $(p,q)$-estimates and
$(E^*_i)$ satisfying $(q',p')$-estimates is easy to establish \cite{Pr}.  
Half of the tree estimate duality a) $\Leftrightarrow$
d)  follows from Lemma~\ref{lem3.5}, which proves that if $X$
satisfies $(\infty,q)$ estimates then $X^*$ satisfies
$(q',1)$-estimates, and if $X^*$ satisfies $(\infty,p')$-tree estimates
$X$ satisfies $(p,1)$-tree estimates. 
But we do not have a direct proof of the other half,
i.e., without first showing (a) $\iff$ (b) and then using Prus'
result, which shows that if $X$ satisfies $(p,1)$-tree estimates
then $X^*$ satisfies $(\infty,p')$ estimates.

Theorem~\ref{thm4.1} was proved in \cite{OS2} and rather than just
repeat that proof we shall give a sketch of the proof emphasizing
the new ideas necessary to go beyond the proof of 
Theorems~\ref{thm:main}  and \ref{thm:equivconditions}.
But first let's see what is an easy consequence of our earlier arguments.

First consider the case where $X$ satisfies $(p,1)$-tree estimates
for some $1<p<\infty$. 
Let $X\subseteq Z$, a reflexive space with a
basis (by 1.3). From Theorem~\ref{thm:main} there exists a blocking
$E=(E_i)$ of the basis for $Z$ so that $X$ naturally embeds into
$Z_p(E)$. $E$ is a bimonotone FDD for $Z_p(E)$ which satisfies
1-$(p,1)$-estimates and thus is boundedly complete. Let $\mathcal{F}
= \{\sum a_if_i\colon \ (a_i)\in B_{\ell_{p'}}$ and $(f_i)$ is a
(finite or infinite) block sequence of $(E^*_n)$ in $S_{Z^*}\}$. It
is easy to check that $\mathcal{F}$ is a $w^*$-compact 1-norming
(for $Z_p(E)$) subset of $B_{Z_p(E)^*}$ and thus $Z_p(E)$ embeds
isometrically into $C(\mathcal{F})$. Furthermore it is again easy to
check that each normalized block sequence of $E$ in $Z_p(E)$ is
pointwise null on $\mathcal{F}$. Hence $E$ is shrinking in $Z_p(E)$
and so $Z_p(E)$ is reflexive. Note for later that this argument only
requires that $E$ is a shrinking FDD.

So we have proved part of Theorem~\ref{thm4.1} in a special case.
Assume now that $X$ satisfies $(p,p)$-tree estimates. In this case
things become simpler. We could follow the proof of 
Theorem~\ref{thm:main} but after obtaining the FDD $E$ for $Z$ so that all
$\bar\delta$-skipped block sequences of $E$ is $S_X$ $2C$-dominate
the unit vector basis of $\ell_p$ we could repeat the argument for
upper estimates and by blocking again obtain an FDD, still denoted
by $E$, so that such $\bar\delta$-skipped block sequences are also
$2C$-dominated by the unit vector basis of $\ell_p$. Then by
estimates as in the proof of Theorem~\ref{thm:main} we could show
that $X$ naturally embeds into $(\sum F_n)_{\ell_p}$ for some
blocking $(F_n)$ of $(E_n)$.

The more general cases of Theorem~\ref{thm4.1} present new difficulties. 
The norm defining $Z_p(E)$ yields $(p,1)$-estimates.
There seems to be no natural way however to directly define a norm
yielding $(\infty,q)$-estimates. 
However if $X$ satisfies $(\infty,q)$-tree estimates then by 
Lemma~\ref{lem3.5} $X^*$ satisfies $(q',1)$-tree estimates. 
We thus need to show that $X^*$
is a quotient of a reflexive space $Y_{q'}(F)$ and obtain $X$ embeds
into $Z= Y_{q'}(F)^*$ which, as is easily seen, satisfies
$(\infty,q)$-estimates for the FDD $(F^*_i)$. 
Then we use the ``$X$
embeds into $Z_p(G^*)$'' argument above, for some blocking $(G^*_i)$
of $(F^*_i)$, to obtain $X$ embeds into a space with an FDD
satisfying $(p,q)$-estimates. 
Of course it needs to be checked that $Z_p(G^*)$ 
preserves the $(\infty,q)$-estimates. 
This was proved by Prus. 
In fact if $F=(F_i)$ is an FDD for $Z$ satisfying
$C$-$(\infty,q)$-estimates then $F$ satisfies
$C$-$(\infty,q)$-estimates for $Z_p(E)$. 
We will not give the proof
but note that the same argument (due to Johnson and Schechtman) is
used below in the proof of Theorem~\ref{T:5.3} (see the Remark after
the proof of Theorem~\ref{T:5.3}).

We thus require the following theorem of which part a) has been
proved.

\begin{thm}\label{thm4.2}
Let $X$ be a reflexive space and let $1<p<\infty$. If $X$ satisfies
$(p,1)$-tree estimates then
\begin{itemize}
\item[a)] If $X$ is a subspace of a reflexive space $Z$ with an FDD $E$ then
there is a blocking $F=(F_i)$ of $E$ so that $X$ naturally embeds
into the reflexive space $Z_p(F)$.
\item[b)] $X$ is a quotient of a reflexive space with an FDD satisfying
$(p,1)$-estimates.
\end{itemize}
\end{thm}

Theorem~\ref{thm4.1} follows readily from Theorem~\ref{thm4.2} (and
Lemma~\ref{lem3.5}). 
We are left with the

\begin{proof}[Sketch of the proof of Theorem~\ref{thm4.2} {\rm b)}]
By Lemma~\ref{L:2.7} we can regard $X^*\subseteq Z^*$ where $Z^*$ is
a reflexive space with a bimonotone FDD $(E^*_i)$ such that
$c_{00}\left(\oplus^\infty_{i=1} E^*_i\right) \cap X^*$ is dense in
$X^*$. Thus we have a quotient map $Q\colon\ Z\to X$. By part a),
$X\subseteq W$, a reflexive space with an FDD $(F_i)$ satisfying
$C$-$(p,1)$-estimates for some $C$.

By a fundamental blocking  lemma of Johnson and Zippin \cite{JZ1} we
may assume that for all $i\le j$, $Q\left(\oplus_{n\in (i,j]}
E_n\right)$ is essentially contained in $\oplus_{n\in [i,j]} F_n$.

We shall increase the norm on $Z$, obtaining a space $\widetilde Z$
for which $(E_i)$, now designated $(\widetilde E_i)$, remains a
shrinking FDD and so that $Q$, now called $\widetilde Q$, remains a
quotient map. Then we shall find a blocking $\widetilde H$ of
$\widetilde E$ so that $\widetilde Q\colon \ \widetilde
Z_p(\widetilde H)\to X$ remains a quotient map. As noted above
$\widetilde Z_p(\widetilde H)$ is reflexive, since $(\widetilde H)$
is shrinking.

For $z\in E_i$ we set $|||\tilde z||| = \|Q(z)\|$ and more generally
for $\tilde z =  \sum \tilde z_i\in c_{00} (\oplus^\infty_1
\widetilde E_i)$ we set $|||\tilde z||| = \max_{m\le n}
\|\sum^n_{i=m} Q(z_i)\|$. Then one checks that
$\widetilde Q$ remains a precise quotient map from $\widetilde Z=$
completion of $c_{00}( \oplus^\infty_1 \widetilde E_i)$
under $|||\cdot|||$ onto $X$. In fact if $Qz=x$, $\|z\| = \|x\|$,
then $|||\tilde z||| = \|z\|$, $\widetilde Q\tilde z = x$. Also
$(\widetilde E_i)$ is a bimonotone FDD for $\tilde z$ (by blocking
we may assume $\widetilde E_i\ne \{0\}$).

A key feature of $(\widetilde Z, |||\cdot|||)$ is the following which is 
easily verified. 
\begin{align}\label{E:4.1}
&\text{If $(\tilde z_i)$ is a block sequence of $(\widetilde E_i)$
in $B_{\widetilde Z}$
  and $(Q\tilde z_i)$ is a basic sequence in $X$ with}\\
&\text{projection constant $\overline K$ and $a\equiv \inf\limits_i
\|\widetilde Q\tilde z_i\|>0$ then}\notag\\
&\qquad\left\|\sum a_i\widetilde Q(\tilde z)_i)\right\| \le
\left|\left|\left| \sum a_i\tilde z_i\right|\right|\right| \le
\frac{3\overline K}a \left\|\sum a_i \widetilde Q(\tilde
z_i)\right\|\notag
\end{align}
for all scalars $(a_i)$.

{From} \eqref{E:4.1} and the fact that
$c_{00}\left(\oplus^\infty_{i=1} E^*_i\right)\cap X^*$ is dense in
$X^*$ one can deduce that $(\tilde E_i)$ is a shrinking FDD for
$\widetilde Z$.

It remains only to prove that there exists $A<\infty$ and a blocking
$\widetilde H$ of $\widetilde E$ satisfying the following. Let $x\in
S_X$. There exists $\tilde z = \sum\tilde z_i$, $\tilde z_i\in
\widetilde H_i$, so that if $(\tilde w_n)$ is any blocking of
$(\tilde z_i)$ then $(\sum |||\tilde w_n|||^p)^{1/p}\le A$ and
$\|\widetilde Q\tilde z-x\|< 1/2$. Thus $\widetilde Q\colon \
\widetilde Z_p(\widetilde H)\to X$ remains a quotient map.

To accomplish this we first use the Johnson and Zippin \cite{JZ1}
 blocking lemma
 for our original $Q\colon \ Z\to X$
to produce a blocking $(C_n)$ of $(E_n)$, and corresponding blocking
$(D_n)$ of $(F_n)$ so that if $x\in S_X$ is essentially contained in
$\oplus_{s\in (i,j)} D_s$ then there exists $z\in B_Z$ with
$Qz\approx x$ and $z\in C_{i,R}\oplus \left(\oplus_{s\in (i,j)}
C_s\right)\oplus C_{j,L}$ where $C_{i,R}$ is the ``right half'' of
the blocking of $E_i$'s yielding $C_i$ and $C_{j,L}$ is the ``left
half'' of $C_j$.

Then we use Lemma~3.2 for suitable $(\delta_i)$ to obtain a blocking
$(G_n)$ of $(D_n)$ and let $(H_n)$ be the corresponding blocking of
$(C_n)$. If $x\in S_X$ we write $x=\sum x_i$, $(x_i)\subseteq X$,
as in Lemma~3.2 and let
\[
B = \{i\colon \ \|P^D_{(t_{i-1},t_i)} x_i-x_i\| <\delta_i\|x_i\|\},
\]
$y=\sum_{i\in B} x_i$. Then $\|y-x\|<1/4$ if $\sum \delta_i<1/4$. 
{From} our left half/right half construction
above we can choose a block sequence $(z_i)_{i\in B}$ of $(E_n)$ in
$B_Z$ with $\|Qz_i-\bar x_i\|\approx 0$ for $i\in B$ and $\bar x_i =
x_i/\|x_i\|$. 
$(\bar x_i)_{i\in B}$ is a perturbation of a block
sequence of $(F_i)$ in $W$ and so admits $2C$-$(p,1)$-estimates. 
{From} \ref{E:4.1} $(\tilde z_i)_{i\in B}$ is equivalent to 
$(\bar x_i)_{i\in B}$ and if we set $\tilde z 
= \sum_{i\in B}\|x_i\|\tilde z_i$ we can show this has the desired
property.
\end{proof}

Suppose that $X$ is a reflexive space which can be renormed to be
a.u.s.\ and can also be renormed to be a.u.c. From 
Theorem~\ref{thm:equivconditions} it follows that there exists $1<q\le
p<\infty$ so that $X$ satisfies $(p,q)$-tree estimates. Thus we have
from Theorem~\ref{thm4.1}.

\begin{thm}\label{thm4.4}
If $X$ is a reflexive space with an equivalent a.u.s.\ norm and an
equivalent a.u.c.\ norm then there exists $1<q\le p<\infty$ so that
\begin{itemize}
\item[a)] $X$ embeds into a reflexive space with an FDD satisfying
1-$(p,q)$-estimates. Hence
\item[b)] $X$ can be renormed to be simultaneously a.u.s.\ of power type $q$ and
a.u.c.\ of power type $p$.
\end{itemize}
\end{thm}

\begin{remark}
The hypothesis of Theorem~\ref{thm4.4} is equivalent to:\ $X$ is
reflexive and $\Sz(X) = \Sz(X^*) = w$.
\end{remark}

It is natural to ask if the results obtained above for
$(p,q)$-estimates can be extended to more general estimates, say
where $\ell_p$ is replaced by a space $V$ with a normalized
1-unconditional basis $(v_i)$ replacing the unit vector basis of
$\ell_p$ and similarly for $\ell_q$. This is done in \cite{OSZ}. The
arguments have a similar flavor as do the ones above but the proofs
are more technically difficult. The analog of Theorem~\ref{thm4.2}
is the following result. The definitions are the analogs of the ones
in the $\ell_p$-case.

\begin{thm}\label{thm4.5}\cite[Theorem 3.1]{OSZ}
Let $V$ be a Banach space with a normalized 1-unconditional basis
$(v_i)$ and let $X$ be a reflexive space satisfying $V$-lower tree
estimates (i.e., for some $C<\infty$ every weakly null tree in $S_X$
admits a branch $C$-dominating $(v_i)$). Then
\begin{itemize}
\item[a)] For every reflexive space $Z$ with an FDD $E=(E_i)$ containing $X$
there is a blocking $H=(H_i)$ of $E$ so that $X$ naturally embeds
into $Z_V(H)$.
\item[b)] There is a space $Y$ with a shrinking FDD $G$ so that $X$ is a
quotient of $Y_V(G)$.
\end{itemize}
\end{thm}

The norm in $Z_V(H)$ is given by for $x\in
c_{00}\left(\oplus^\infty_1 H_i\right)$ by
\[
\|x\|  = \sup\bigg\{\Big\|\sum^\infty_{j=1}\|P^H_{(n_{j-1},n_j]}
x\|_Z \cdot v_i\Big\|_V\colon\ 0 = n_0<n_1<\cdots\bigg\}.
\]
Unlike the $\ell_p$ case $(H_i)$, which is an FDD for $Z_V(H)$, does
not automatically admit a lower $V$-estimate on blocks. But this can
be achieved with additional hypotheses on $V$.

\begin{Defn}
A normalized 1-unconditional basis $(v_i)$ is \emph{regular} iff
\begin{itemize}
\item[i)] $(v_i)$ is dominated by every normalized block basis of $(v_i)$.
\item[ii)] There exists $c>0$ so that for all $(a_i)\in c_{00}$ and $n\in
\mathbb{N}$
\[
\Big\|\sum^\infty_{i=1} a_iv_{i+n}\Big\| \ge
c\Big\|\sum^\infty_{i=1} a_iv_i\Big\|.
\]
\item[iii)] There exists $d>0$ so that for all $m\in \mathbb{N}$ there exists $L
= L(m)\ge m$ so that for all $k\le m$
\[
\Big\|\sum^\infty_{i=L+1} a_iv_{i-k}\Big\| \ge
d\Big\|\sum^\infty_{i=L+1} a_iv_i\Big\| \quad\text{whenever}\quad
(a_i)\in c_{00}.
\]
\end{itemize}
\end{Defn}

\begin{thm}\label{thm4.6}\cite[Corollary 3.2]{OSZ}
Let $V$ be a reflexive space with a regular normalized
1-uncondi\-tional basis $(v_i)$. Let $X$ be a reflexive space with
$V$-lower tree estimates. Then $X$ is a subspace of a reflexive
space $Z$ with an FDD satisfying $V$-lower estimates and $X$ is a
quotient of a reflexive space $Y$ with an FDD satisfying $V$-lower
estimates.
\end{thm}

For an upper and lower estimate result we have

\begin{thm}\label{thm4.7}\cite[Theorem 3.4]{OSZ}
Let $V$ and $U^*$ be reflexive Banach spaces with regular normalized
1-unconditional bases $(v_i)$ and $(u^*_i)$, respectively. Assume
that every subsequence of $(u_i)$ dominates every normalized block
basis of $(v_i)$ and every normalized block basis of $(u_i)$
dominates every subsequence of $(v_i)$. If $X$ is a reflexive space
satisfying $(V,U)$-tree estimates then $X$ embeds into a reflexive
space $Z$ with an FDD satisfying $(V,U)$-estimates.
\end{thm}

Examples of spaces $(V,U)$ satisfying the hypothesis of 
Theorem~\ref{thm4.7} are the {\em convexified Tsirelson} spaces
$(T_{p,\gamma},T^*_{q',\gamma})$ where $1\le q\le p\le \infty$ and
$0<\gamma<1/4$. If $X$ is a reflexive asymptotic $\ell_p$ space
(i.e., $\exists\  C\ge 1\ \forall\ n$ $\forall\ (e_i)^n_1\in \{X\}_n$,
\[
C^{-1}\bigg(\sum^n_1|a_i|^p\bigg)^{1/p} \le \Big\|\sum^n_1
a_ie_i\Big\| \le C\bigg(\sum^n_{i=1} |a_i|^p\bigg)^{1/p}
\]
for all $(a_i)^n_1\subseteq \mathbb{R}$) then it can be easily seen
that $X$ satisfies $(T_{p,\gamma},T^*_{p'\gamma})$-tree estimates
for some $0<\gamma < 1/4$. As an application we have

\begin{cor}\label{cor4.8}
Let $X$ be a reflexive asymptotic $\ell_p$ space. Then $X$ embeds
into a reflexive space with an asymptotic $\ell_p$ FDD. $X$ is also
a quotient of such a space.
\end{cor}

Similar results can  be obtained analogous to those of 
Theorem~\ref{thm4.1} (see \cite{OSZ}).

\section{Universal spaces}\label{S:5}

We begin with the solution to Bourgain's problem (see Section
\ref{S:1}). Note that (e.g., by Krivine's theorem \cite{K}) if $X$
contains an isomorph of $\ell_p$ for all $1<p<\infty$, then $c_0$
and $\ell_1$ are finitely represented in $X$ so $X$ cannot be
superreflexive.

One step in the proof will be, given $1<q\le p<\infty$, to construct
a space
 $Z_{(p,q)}$ with an FDD satisfying $1$-$(p,q)$-estimates which is
 universal for all such spaces. We shall do this first before
 proceeding to the theorem.
 S.~Prus \cite{Pr} has shown  a similar result but we prefer to
  present a somewhat different argument which could prove useful
   elsewhere.

  \begin{lem}\label{L:5.1}
  Let $1\le q\le p\le\infty$, and
   let $F$ and $G$ be two finite dimensional normed linear
   spaces. Denote the norm on $F$ and $G$ by
     $\|\cdot\|_E$ and $\|\cdot\|_G$ respectively.
 Let $|||\cdot|||$ be a norm on $F\oplus G$ and assume that
 $(F,G)$ satisfies $1$-$(p,q)$-estimates in $(F\oplus G,|||\cdot|||)$
 and
  there are
$1<c<d<\infty$ so that
\begin{align}\label{E:5.1.1}
&c\|f\|_F\le |||f|||\le d\|f\|_F\text{ whenever }f\in F\text{ and }\\
&c\|g\|_G\le |||g|||\le d\|g\|_G\text{ whenever }g\in G.\notag
\end{align}
 Then there is a norm $\|\cdot\|$ on $F\oplus G$ extending
 $\|\cdot\|_F$ and $\|\cdot\|_G$, so that $(F,G)$ is an
  FDD  satisfying $1$-$(p,q)$-estimates in
  $(F\oplus G,\|\cdot\|)$ and
\begin{equation}\label{E:5.1.2}
 c\|f+g\|\le |||f+g|||\le d\|f+g\|\text{ whenever }f\in F\text{ and } g\in G.
\end{equation}
\end{lem}

\begin{proof} For $f\in F$ and $g\in G$ put:
\begin{equation*}
\|f+g\|=\max\Big\{
\big(\|f\|_F^p+\|g\|_G^p\big)^{1/p},\frac1d|||f+g|||\Big\}\ ,
\end{equation*}
where we replace $(\|f\|_F^p+\|g\|_G^p\big)^{1/p}$ by
$\max\{\|f\|_F,\|g\|_G\}$ if $p=\infty$. Clearly $(F,G)$ satisfies
$1$-$(p,1)$-estimates, and since $(F,G)$ satisfies
$1$-$(\infty,q)$-estimates in $(F\oplus G,|||\cdot|||)$,
 this is also true for  $(F\oplus G,\|\cdot\|)$.
Moreover, for $f\in F$ and $g\in G$ we deduce that 
\begin{align*}
c\|f+g\|&=\max\Big\{ \big(\|cf\|_F^p+\|cg\|_G^p\big)^{1/p},
\frac{c}{d}|||f+g|||\Big\}\\
 &\le \max\Big\{\big(|||f|||_F^p+|||g|||_G^p\big)^{1/p},|||f+g|||\Big\}\\
 &\le |||f+g|||\le d\|f+g\|\ .
\end{align*}
\end{proof}
We introduce the following terminology.
\begin{defn}
 Let $E_\alpha$ be finite
 dimensional linear space for each $\alpha\in T_\infty$ and
  let $\|\cdot\|_\beta$ be a norm
   on $\coo(\oplus_{i=1}^\infty E_{\beta_i})$ for each
    branch $\beta=(\beta_i)_{i=1}^\infty$ of $T_\infty$.
     We say that the family $(\|\cdot\|_\beta)$ indexed over
      all branches of $T_\infty$ is {\em compatible} if
\begin{enumerate}
\item For every branch $\beta=(\beta_i)_{i=1}^\infty$ of $T_\infty$,
 $(E_{\beta_i})$ is a bimonotone FDD for the completion $X_\beta$
  of $\coo(\oplus_{i=1}^\infty E_{\beta_i})$ under $\|\cdot\|_\beta$.
\item  If $\alpha=(\alpha_i)$ and $\beta=(\beta_i)$ are two branches
 and if $\ell=\max\{i: \forall\  j\le i\,\,\alpha_j=\beta_j\}$ ($\ell=0$ if
 $\alpha_1\not=\beta_1$) then
 $\|\cdot\|_\alpha$ and $\|\cdot\|_\beta$ coincide on
  $\oplus_{i=1}^\ell E_{\alpha_i}$.
\end{enumerate}
\end{defn}
\begin{prop}\label{P:5.2}
Let $1\le q\le p\le \infty$. Then there exists a tree
$(E_\alpha)_{\alpha\in T_\infty}$ of finite dimensional linear
spaces and a compatible family of norms $\|\cdot\|_\beta$ for each
 branch $\beta$ of $T_\infty$ satisfying the following
\begin{enumerate}
\item If $\beta=(\beta_i)$ is a branch in $T_\infty$
 then $(E_{\beta_i})$ satisfies $1$-$(p,q)$-estimates for
 $\|\cdot\|_\beta$.
\item Let $Y$ be any Banach space with norm $\|\cdot\|$ and
 with an FDD $(F_i)$ satisfying $1$-$(p,q)$-estimates in $Y$ and
  let $d>1$. Then there exists a branch $\beta=(\beta_i)$ of
  $T_\infty$ and an isomorphism $I$ from
   $X_\beta$ (the completion of $\coo(\oplus E_{\beta_i})$) onto $Y$ under
    $\|\cdot\|_\beta$) mapping $E_{\beta_i}$ onto $F_i$, for $i\in\nat$, satisfying
$$\|x\|_\beta\le \|I(x)\|\le d \|x\|_\beta\text{ whenever }
 x\in X_\beta.$$
\end{enumerate}
\end{prop}
\begin{proof} For $n\in\nat$ let $T_n$ be the elements of $T_\infty$
of length $n$.
 By induction on $n\in\nat$ we will
 define the normed linear spaces $E_\alpha$
  for all $\alpha\in T_n$ and
   norms $\|\cdot\|_\beta$ on $\oplus_{j=1}^n E_{(\alpha_j)}$
 where $\beta=(\alpha_1,\alpha_2,\ldots,\alpha_n)$ is a branch of
  length $n$ in $\bigcup_{j=1}^n T_j$, i.e. for $i\le n$  $|\alpha_i|=i$
   and $\alpha_i$ is a successor of $\alpha_{i-1}$ if $1<i$.

The first level of $(E_\alpha)_{\alpha\in T_\infty}$ is any sequence
of finite dimensional Banach spaces which is dense (with respect to
the Banach-Mazur distance) in the set of all finite dimensional
Banach spaces.

Assume we have defined for
  a branch
 $\beta=(\alpha_1,\alpha_2,\ldots,\alpha_n)$ the space
 $E^\beta=\oplus_{i=1}^n E_{\alpha_i}$ along with a norm
 $\|\cdot\|_\beta$ on it. Let
$\beta=(\alpha_1,\alpha_2,\ldots, \alpha_n)$  be such a branch. The
successors of $\alpha_n$ are chosen as follows. Let $(G_i)$ be the
 spaces of level 1. For each $G_i$
 we consider the set of all extensions of $\|\cdot\|_\beta$ to
  $E^\beta\oplus G_i$ satisfying $1$-$(p,q)$-estimates.

For any two such extensions $\|\cdot\|_1$ and $\|\cdot\|_2$ we
define
 the distance between $\|\cdot\|_1$ and $\|\cdot\|_2$ by
 \begin{equation*}
 d(\|\cdot\|_1,\|\cdot\|_2)=\ln\big(\|I\|\cdot\|I^{-1}\|\big),
 \end{equation*}
 where $I:\big(E_\alpha\oplus G_i,\|\cdot\|_1\big)\to
 \big(E_\alpha\oplus G_i,\|\cdot\|_2\big)$ is the identity.
 We then
 choose a countable dense subset of these extensions
 with respect to $d(\cdot,\cdot)$.
 The sequence of all successors will then be formed
   by the union over all $i$ of these countable many
  extensions.

  To see (2) we will use Lemma~\ref{L:5.1}. Let $1<d$ and let
   $1<c_n<d_n<d$ with $c_n\searrow 1$ and $d_n\nearrow d$, if
   $n\nearrow\infty$. Let $Y$ and $(F_i)$ as in (2 and denote
 the norm on $Y$ by $\|\cdot\|$. To start
    we find $\alpha_1\in T_1$ and an isometry  $I_1$
 from $F_1$ onto $(E_{\alpha_1},|||\cdot|||_1)$
 where $|||\cdot|||_1$ is a norm on $E_{\alpha_1}$ with
  $c_1\|x\|_{\beta_1}\le |||x|||_1\le d_1\|x\|_{\beta_1}$
   ($\beta_1=(\alpha_1)$).

Assume we constructed a branch $\beta=(\alpha_1,\ldots, \alpha_n)$
 of length $n$ along with a norm $|||\cdot|||_n$ on $E^\beta$ and an
 isometry mapping, $F_i$ onto $E_{\alpha_i}$, for $i=1,2,\ldots, n$
 $$I_n:\big(\oplus_{i=1}^n F_i,||\cdot||\big)\to
 \big(E^\beta,|||\cdot|||_n\big)$$
 satisfying
 $$c_n\|x\|_\beta\le|||x|||_n\le d_n\|x\|_\beta
  \text{ for }x\in E^{\beta}\ .$$
Since $(G_i)$ is dense in the set of all finite dimensional normed
spaces
 we can find a $G=G_i$, whose norm we denote by $\|\cdot\|_G$,  $\dim(G)=\dim(F_{n+1})$ and an isometry
 $J:F_{n+1}\to (G,|||\cdot|||)$ where $|||\cdot|||$ is a  norm on $G$ satisfying
   $$ c_n\|x\|\le |||x|||\le d_n\|x\|_G,\text{ whenever }
    x\in G\ .$$
Define
$$I_{n+1}: \oplus_{i+1}^n F_i\to E^\beta\oplus G,
\quad \sum_{i=1}^{n+1} x_i\mapsto I_n\Big(\sum_{i=1}^{n}x_i\Big)+J(x_{n+1})\ ,$$
and put $$|||x+y|||_{n+1}=||I_{n+1}^{-1}(x,y)||\text{ whenever $x\in
E^\beta$ and $y\in G$}.$$

By Lemma~\ref{L:5.1} we can find a norm $\|\cdot\|$ on
$E^\beta\oplus G$
 extending $\|\cdot\|_\beta$ on $E^\beta$ and $||\cdot||_G$
  on $G$  for which $(E^\beta,G)$ satisfies $1$-$(p,q)$-estimates
   so that
 $$ c_n\|x+y\|\le |||x+y|||_{n+1}      \le d_n\|x+y\|
\text{ whenever $x\in E^\beta$ and $y\in G$.}$$ {From} our
construction of $(E_\alpha)_{\alpha\in T_\infty}$
 there exists a  successor $\alpha_{n+1}$ of $\alpha$ so that for
$\overline{\beta}=(\alpha_1,\alpha_2,\ldots,\alpha_{n+1})$ and
$x\in \oplus E^{\overline{\beta}}$
$$c_{n+1}\|x\|_{\overline{\beta}}\le |||(x)|||_{n+1}
\le d_{n+1}\|x\|_{\overline{\beta}}\ ,$$
which finishes our recursive choice.

Taking now the infinite branch $\beta=(\alpha_i)_{i=1}^{\infty}$
yields our
 claim (2).
\end{proof}
\begin{thm}\label{T:5.3}
There exists a separable reflexive space $X_u$ which is universal
for $\{X:X$ is reflexive and $\Sz(X)=\Sz (X^*)=\w\}$. In particular
$X_u$ contains an isomorph of all separable superreflexive spaces.
\end{thm}

\begin{proof}
We first note that if $X$ is superreflexive then $X$ satisfies
$(p,q)$-tree estimates for some $1<q\le p<\infty$ (\cite{Ja},
\cite{GG}). By Theorem~\ref{thm4.1} $X$ then embeds into a reflexive
space $Z$ with an FDD satisfying $1$-$(p,q)$-estimates. Moreover by
Theorem~\ref{thm:equivconditions}  (applied to  $X$ and $X^*$) the
same holds if $X$ is reflexive with $\Sz(X) = \Sz (X^*) =\w$. Thus
it suffices to produce a space $Z_{(p,q)}$ with an FDD satisfying
1-$(p,q)$-estimates which is universal for all spaces with a
1-$(p,q)$ FDD. We then take $X_u = (\sum Z_{(p_n,q_n)})_{\ell_2}$
where $p_n \uparrow \infty$ and $q_n\downarrow 1$.

 To construct $Z_{(p,q)}$ we first let
 $(E_\alpha)_{\alpha\in T_\infty}$ along with compatible norms
 $\|\cdot\|_\beta$ for branches $\beta$ in $T_\infty$ be as
 constructed in Proposition~\ref{P:5.2} for $p$ and $q$.

$Z_{(p,q)}$ is then the completion of $c_{00}(\oplus_{\alpha\in
T_\infty} E_\alpha)$ under
$$\|z\| = \sup \bigg\{ \Big( \sum_j \|P_{I_j}^E z\|^p\Big)^{1/p}
: I_1, I_2,\ldots, \text{ are disjoint segments in }
T_\infty\bigg\}\ .$$ For a segment $I$, $\|P_I^E z\| = \| P_I^E
z\|_\beta$ where $\alpha$ is any branch containing $I$. $E=
(E_\alpha)_{\alpha\in T_\infty}$ is thus an FDD for $Z_{(p,q)}$ when
ordered linearly in any manner compatible with the tree order of
$T_\alpha$, e.g., $E_{(n_1,n_2,n_3)}$ comes after $E_{(n_1,n_2)}$.
Moreover $E$, when thus ordered, satisfies 1-$(p,1)$-estimates and
the norm on each branch of $(E_\alpha)_{\alpha\in T_\alpha}$ is
preserved. Finally we check that $Z_{(p,q)}$ satisfies
1-$(\infty,q)$-estimates.

Let $z = \sum z_i \in c_{00}(\oplus_{\alpha\in T_\infty} E_\alpha)$
 where $(z_i)$ is a block sequence of $E$. Let $\|z\| = (\sum\|
P_{I_j}^E z\|^p )^{1/p}$ where $I_1,I_2,\ldots$ are disjoint
segments in $T_\infty$. We decompose $I_j$ into segments
$I_{j,1},I_{j,2},\ldots$ so that $P_{I_{j,i}}^E z_i = P_{I_{j,i}}^E
z$ and $P_{I_{j,i}}^E z_s = 0$ if $s\ne i$. Then
$$\|z\| \le \bigg( \sum_j \Big( \sum_i \| P_{I_{j,i}} z_i\|^q\Big)^{p/q}
\bigg)^{1/p}$$ by the 1-$(\infty,q)$-estimates on each branch. Now
$$\bigg( \sum_j \Big(\sum_i \|P_{I_{j,i}} z_i\|^q\Big)^{p/q}\bigg)^{1/p}
\le \bigg[ \sum_i \Big(\sum_j \|P_{I_{j,i}}
z_i\|^p\Big)^{q/p}\bigg]^{1/q}$$ by the reverse triangle inequality
in $\ell_{p/q}$. Thus $\|z\| \le (\sum_i \|z_i\|^q)^{1/q}$.
\end{proof}

\begin{remark} 
The clever argument for the $1$-$(\infty,q)$ estimate is due
to Johnson and Schechtman. It was used in \cite{OS2} to show that if
 an FDD $E=(E_i)$ for a space $Z$ satisfies $1$-$(\infty,q)$-estimates
  then it also satisfies $1$-$(\infty,q)$-estimates in $Z_p(E)$.
\end{remark}

 We now turn to the universal problem for the classes (see
  Theorem~\ref{thm:equivconditions} )
\begin{align*}\C_{\auc} &= \{ Y : Y \text{ is separable, reflexive and has an
equivalent a.u.c. norm}\}\\
&=\{ Y:Y \text{ is separable, reflexive and $\Sz(Y^*)=\omega$}\}\\
\intertext{and} \C_{\aus}&= \{ Y : Y \text{ is separable, reflexive
            and has an equivalent a.u.s. norm}\}\\
  &=\{ Y:Y \text{ is separable, reflexive and
                              $\Sz(Y)=\omega$}\}\ .
\end{align*}
First note that the Tsirelson space
$T(\frac12,S_\alpha)$, $\alpha <\w_1$, $S_\alpha = \alpha^{\text{\rm th}}$
Schreier class \cite{AGR}  are all in $\C_{\auc}$ since their unit
vector basis has $(p,1)$-estimates for all $p>1$ and their duals are
all in $\C_{\aus}$. It follows by index arguments
  \cite{Bo}) that any space universal for $C_{\auc}$ must contain
$\ell_1$ and any space universal for $\C_{\aus}$ must contain $c_0$.

\begin{prop}\label{P:5.4}
There exists a separable dual space $X$ which is universal for
$\C_{\auc}$. $X$ is the $\ell_2$ sum of a.u.c. spaces.
\end{prop}

\begin{proof}
The argument is much the same as that of Theorem~\ref{T:5.3}.
 For $p<\infty$ we let $Z_{(p,1)}$ be the space constructed in the
 proof of Theorem~\ref{T:5.3}. $Z_{(p,1)}$ has an FDD satisfying
 $1$-$(p,1)$-estimates and as such, having a boundedly complete FDD,
  is a separable dual space. By Theorem~\ref{thm4.1} $Z_{(p,1)}$ is
   universal for all spaces in $\C_{\auc}$ satisfying
   $(p,1)$-estimates. Thus by Theorem~\ref{thm:equivconditions},
    $X=\big(\oplus_{n=2}^\infty Z_{(n,1)} \big)_{\ell_2}$ is universal for
     $\C_{\auc}$. $X$ is a separable dual space.
\end{proof}

\begin{remark}
The spaces $Z_{(p,q))}$ constructed in Theorem~\ref{T:5.3} and
Proposition~\ref{P:5.4} are actually complementably universal for
the members of their respective classes which have $(p,q)$ or
$(p,1)$ FDD's.
 The space $X$ of Proposition~\ref{P:5.4} is universal for the class
 \begin{align*}
  &\{ Y: Y \text{ is a separable dual
 satisfying $w^*$-$(p,1)$-estimates for some }p<\infty\}\\
 &\quad=\{ Y: Y=W^*\text{ with } \Sz(W)=\omega\}\ .
\end{align*}
\end{remark}

\begin{prop}\label{P:5.5}
There exists a space $Y$ with separable dual which is universal for
the class $\C_{\aus}$. $Y$ is the $\ell_2$ sum of a.u.s. spaces.
\end{prop}

\begin{proof} Let $q>1$. Let $(E_\alpha)_{\alpha\in T_\infty}$
 and a compatible set of norms $\|\cdot\|_\beta$ for each branch
  $\beta$ of $T_\infty$ be constructed as in 
  Proposition~\ref{P:5.2} for $(\infty,q)$. We let $Z_q$ be the completion of
    $\coo(\oplus_{\alpha\in T_\infty} E_\alpha)$ under the norm
     $$\|z\|=\sup\big\{\|P^E_\beta z\|:\beta\text{ is a branch in }
     T_\infty\big\}\ .$$
If $(E_\alpha)_{\alpha\in T_\infty}$ is linearly ordered in a manner
 compatible with the order on $T_\infty$ it becomes a bimonotone
  FDD for $Z_q$ satisfying $1-(\infty,q)$-estimates.

  Let $q_n\searrow 1$, if $n\nearrow\infty$ and set
  $Y=\big(\oplus_{n=1}^\infty Z_{q_n}\big)_{\ell_2}$. By
 Theorem~\ref{thm:equivconditions}
  and Theorem~\ref{thm4.1} $Y$ is universal for $\C_{\aus}$. Clearly
  $Y^*$ is separable.
\end{proof}
\begin{remark} For $Y$ as constructed in Proposition~\ref{P:5.5} it
follows that $\Sz(Y)=\omega^2$.
\end{remark}

Finally we note the following result from \cite{OSZ}.

\begin{thm}\label{thm:5.7}
Let $K<\infty$ and $1\le p\le \infty$. 
There exists a reflexive asymptotic $\ell_p$ space which is universal 
for the class of all reflexive $K$-asymptotic $\ell_p$ spaces.
\end{thm}

We refer to \cite{OSZ} for the proof and for more general versions 
of this results.

\end{document}